\documentclass[a4paper]{amsart}

\usepackage{amssymb,amsmath,graphics,verbatim}
\usepackage{latexsym}
\usepackage{eucal}
\makeatletter
\@addtoreset{equation}{section}\makeatother

\newtheorem{theo}{Theorem}[section]
\newtheorem{lem}[theo]{Lemma}
\newtheorem{prop}[theo]{Proposition}
\newtheorem{cor}[theo]{Corollary}
\newtheorem{defi}[theo]{Definition}

\def\tr{{\rm tr}}
\newcommand{\mc}{\mathcal}

\newcommand{\nn}{\mathbb{N}}
\newcommand{\cc}{\mathbb{C}}
\newcommand{\hh}{\mathbb{H}}
\newcommand{\zz}{\mathbb{Z}}

\newcommand{\la}{\lambda}
\newcommand{\Scal}{{\rm Scal}}
\newcommand{\Ric}{{\rm Ric}}
\newcommand{\eps}{\epsilon}

\newcommand{\pl}{\partial}
\newcommand{\x}{\times}

\newcommand{\til}{\widetilde}

\newcommand{\cjd}{\rangle}
\newcommand{\cjg}{\langle}

\newcommand{\demi}{\frac{1}{2}}
\newcommand{\ndemi}{\frac{n}{2}}
\newcommand{\tra}{\textrm{Tr}}

\def\qed{\hfill$\square$}

\newcommand{\dis}{\displaystyle}

\begin{document}
\title[Conformal harmonic forms and Branson-Gover operators]{Conformal harmonic forms, Branson-Gover operators and Dirichlet problem at infinity.}
\author{Erwann Aubry}
\author{Colin Guillarmou}
\address{Lab. Dieudonn\'e\\
Univ. de Nice Sophia-Antipolis \\
Parc Valrose\\
06108 Nice\\ 
FRANCE}
\email{eaubry@math.unice.fr}
\email{cguillar@math.unice.fr}
\begin{abstract}
For odd dimensional Poincar\'e-Einstein manifolds $(X^{n+1},g)$, we study the set of harmonic $k$-forms (for $k<\ndemi$) 
which are $C^m$ (with $m\in\nn$) on the conformal compactification $\bar{X}$ of $X$.
This is infinite dimensional for small $m$ but it becomes finite dimensional if $m$ is large enough, 
and in one-to-one correspondence with the direct sum of the relative cohomology $H^k(\bar{X},\pl\bar{X})$ 
and the kernel of the Branson-Gover \cite{BG} differential operators $(L_k,G_k)$  on the conformal infinity 
$(\pl\bar{X},[h_0])$. In a second time we relate the set of $C^{n-2k+1}(\Lambda^k(\bar{X}))$ forms
in the kernel of $d+\delta_g$ to the conformal harmonics on the boundary in the sense of \cite{BG}, 
providing some sort of long exact sequence adapted to this setting.
This study also provides another construction of Branson-Gover differential operators, 
including a parallel construction of the generalization of $Q$ curvature for forms. 
\end{abstract}
%
\maketitle

\section{Introduction}

Let $(M,[h_0])$ be an n-dimensional compact manifold equipped with a conformal class $[h_0]$.
The $k$-th cohomology group $H^k(M)$ can be identified with $\ker(d+\delta_{h})$ for any $h\in[h_0]$
by usual Hodge-De Rham Theory. However, the choice of harmonic representatives in $H^k(M)$
is not conformally invariant with respect to $[h_0]$, except when $n$ is even and $k=\ndemi$.  
Recently, Branson and Gover \cite{BG} defined new complexes, new conformally invariant spaces of forms and 
new operators to somehow generalize this $k=\ndemi$ case. More precisely, they introduce
conformally covariant differential operators $L_k^{{\rm BG},\ell}$ of order $2\ell$ on the bundle $\Lambda^{k}(M)$ 
of $k$-forms, for $\ell\in\nn$ (resp. $\ell\in\{1,\dots,\ndemi\}$) if $n$ is odd (resp. $n$ is even).
A particularly interesting case is the critical one in even dimension, this is 
\begin{equation}\label{Lk}
L^{\rm BG}_k:=L_k^{{\rm BG},\frac{n}{2}-k}.\end{equation}   
The main features of this operator are that it factorizes under the form $L^{\rm BG}_k=G^{\rm BG}_{k+1}d$ for some operator 
\begin{equation}\label{Gk}
G^{\rm BG}_{k+1}:C^\infty(M,\Lambda^{k+1}(M))\to C^{\infty}(M,\Lambda^k(M))
\end{equation} 
and that $G^{\rm BG}_k$ factorizes under the form $G^{\rm BG}_k=\delta_{h_0}Q^{\rm BG}_k$
for some differential operator 
\begin{equation}\label{Qk}
Q^{\rm BG}_k: C^\infty(M,\Lambda^{k}(M))\cap \ker d\to C^{\infty}(M,\Lambda^{k}(M)) 
\end{equation} 
where $\delta_{h_0}$ is the adjoint of $d$ with respect to $h_0$.  This gives rise to an elliptic complex 
\[\dots \xrightarrow{d} \Lambda^{k-1}(M)\xrightarrow{d} \Lambda^{k}(M)\xrightarrow{L^{\rm BG}_k}\Lambda^{k}(M)
\xrightarrow{\delta_{h_0}}\Lambda^{k-1}(M)\xrightarrow{\delta_{h_0}}\dots \] 
named the \emph{detour complex}, whose cohomology is conformally invariant. 
Moreover, the pairs $(L^{\rm BG}_k,G^{\rm BG}_k)$ and $(d,G^{\rm BG}_k)$ on $\Lambda^k(M)\oplus \Lambda^{k}(M)$ 
are graded injectively elliptic in the sense that $\delta_{h_0} d+dG^{\rm BG}_k$ and 
$L^{\rm BG}_k+dG^{\rm BG}_k$ are elliptic. Their finite dimensional kernel 
\begin{equation}\label{hkhkl}
\mc{H}_L^k(M):=\ker (L^{\rm BG}_k,G^{\rm BG}_k), \quad \mc{H}^k(M):=\ker(d,G^{\rm BG}_k)
\end{equation}
are conformally invariant, the elements of $\mc{H}^k(M)$ are named \emph{conformal harmonics}, providing a type of Hodge theory for conformal structure. 
The operator $Q^{\rm BG}_k$ above generalizes Branson $Q$-curvature in the sense that it satisfies, as operators
on closed $k$-forms, 
\[\hat{Q}^{\rm BG}_k=e^{\mu(2k-n)}(Q^{\rm BG}_k+L^{\rm BG}_k\mu)\]
if $\hat{h}_0=e^{2\mu}h_0$ is another conformal representative.\\ 

The general approach of Fefferman-Graham \cite{FG} for dealing with conformal invariants 
is related to Poincar\'e-Einstein manifolds, roughly speaking it provides a correspondence
between Riemannian invariants in the bulk $(X,g)$ and conformal invariants on the conformal infinity $(\pl\bar{X},[h_0])$ of $(X,g)$, inspired by the identification of the conformal group of the sphere $S^n$ with the isometry group of the hyperbolic
space $\hh^{n+1}$. A smooth Riemannian manifold $(X,g)$ 
is said to be a \emph{Poincar\'e-Einstein manifold} with conformal infinity 
$(M,[h_0])$ if the space $X$ compactifies smoothly in $\bar{X}$ with boundary $\pl\bar{X}=M$, 
and if there is a boundary defining function of $\bar{X}$ and some collar neighbourhood
$(0,\eps)_x\x\pl\bar{X}$ of the boundary such that
\begin{equation}\label{modelform}
g=\frac{dx^2+h_x}{x^2} 
\end{equation}
\begin{equation}\label{ahe}
\textrm{Ric}(g)=-ng+O(x^{\infty})
\end{equation}
where $h_x$ is a one-parameter family of smooth metrics on $\pl\bar{X}$ 
such that there exist some family of smooth tensors $h_x^j$ ($j\in\nn_0$) on $\pl\bar{X}$, depending smoothly 
on $x\in[0,\eps)$ with 
\begin{equation}\label{hx}
\left\{
\begin{array}{l}
 h_x\sim\sum_{j=0}^\infty h_x^j(x^n\log x)^j \textrm{ as }x\to 0 \textrm{ if }n+1 \textrm{ is odd}\\
 h_x \textrm{ is smooth in }x\in [0,\eps) \textrm{  if }n+1 \textrm{ is even}
\end{array}\right.
\end{equation}
\begin{equation}\label{hxh0}
h_x|_{x=0}\in[h_0].
\end{equation} 
The tensor $h^1_0$ is called \emph{obstruction tensor} of $h_0$, it is defined in \cite{FG} and studied further in 
\cite{GRH}.
We shall say that $(X,g)$ is a smooth Poincar\'e-Einstein manifold if $x^2g$ extends smoothly on $\bar{X}$, i.e.
either if $n+1$ is even or $n+1$ is odd and $h_x^j=0$ for all $j>0$. It is proved in \cite{FG2}
that $h^1_0=0$ implies that $(X,g)$ is a smooth Poincar\'e-Einstein manifold.\\ 

The boundary $\pl\bar{X}=\{x=0\}$ inherits naturally from $g$ the conformal class $[h_0]$ of $h_x|_{x=0}$ since
the boundary defining function $x$ satisfying such conditions are not unique. 
A fundamental result of Fefferman-Graham \cite{FG}, which we do not state in full generality, 
is that for any $(M,[h_0])$ compact that can be realized as the boundary of
smooth compact manifold with boundary $\bar{X}$, there is a Poincar\'e-Einstein manifold $(X,g)$ for $(M,[h_0])$, 
and $h_x$ in \eqref{hx} is uniquely determined by $h_0$ up to order $O(x^{n})$ and up to diffeomorphism which restricts
to the Identity on $M$.
The most basic exemple is the hyperbolic space 
$\hh^{n+1}$ which is a smooth Poincar\'e-Einstein for the canonical conformal structure of the sphere $S^n$,
as well as quotients of $\hh^{n+1}$ by convex co-compact groups of isometries.\\ 

It has been proved by Mazzeo \cite{Ma} that\footnote{The class of manifold considered by Mazzeo is actually 
larger and does not require the asymptotic Einstein condition \eqref{ahe}} 
for a Poincar\'e-Einstein manifold $(X,g)$, the relative cohomology $H^k(\bar{X},\pl\bar{X})$
is canonically isomorphic to the $L^2$ kernel $\ker_{L^2}(\Delta_k)$
of the Laplacian $\Delta_k=(d+\delta_g)^2$ with respect to the metric $g$, 
acting on the bundle $\Lambda^k(\bar{X})$ of $k$-forms if $k<\ndemi$. In other terms
the relative cohomology has a basis of $L^2$ harmonic representatives.
In this work, we give an interpretation of the spaces $\mc{H}^k,\mc{H}_L^k$ in terms of harmonic 
forms on the bulk $X$ with a certain regularity on the compactification $\bar{X}$. 

\begin{theo}\label{th1}
Let $(X^{n+1},g)$ be an odd dimensional Poincar\'e-Einstein manifold with conformal infinity $(M,[h_0])$ 
and let $\Delta_k=(d+\delta_g)^2$ be the induced Laplacian on $k$-forms on $X$ where $0\leq k<\ndemi-1$. 
For $m\in\nn$ and $0<k<\ndemi-1$, define 
\[K^k_m(\bar{X}):=\{\omega\in C^{m}(\bar{X};\Lambda^k(\bar{X})); \Delta_k\omega=0\},\]
then $K_m^k(\bar{X})$ is infinite dimensional for $m<n-2k+1$ while it is finite dimensional
for $m\in [n-2k+1,n-1]$ and there is a canonical short exact sequence
\begin{equation}\label{shortexact}
 0 \longrightarrow{} H^{k}(\bar{X},\pl\bar{X})\xrightarrow{i} K_{m}^k(\bar{X})\xrightarrow{r} \mc{H}^k_L(M)\longrightarrow{} 0
\end{equation} 
where $\mc{H}_L^k$ is defined in \eqref{hkhkl} and $H^k(\bar{X},\pl\bar{X})$ is the relative cohomology
space of degree $k$ of $\bar{X}$, $i$ denotes inclusion and $r$ denotes pull back by the natural inclusion
$\pl\bar{X}\to \bar{X}$. 
If in addition the Fefferman-Graham obstruction tensor of $(M,[h_0])$ vanishes, i.e. if $(X,g)$ is a smooth Poincar\'e-Einstein manifold, then $K^{k}_{n-2k+1}(\bar{X})=K_\infty^k(\bar{X})$.\\
When $k=\ndemi-1$, the same results hold by replacing $K_{n-2k+1}^k(\bar{X})$
by the set of  harmonic forms in $C^{n-2k+1,\alpha}(\bar{X},\Lambda^{k}(\bar{X}))$ for some $\alpha\in(0,1)$.\\   
When $k=0$, $K^0_m(\bar{X})$ is infinite dimensional for $m<n$ while $K_n^0(\bar{X})$ is finite dimensional
and the exact sequence \eqref{shortexact} holds.
\end{theo} 

In that purpose, we show that we can recover the Branson-Gover operators $L^{\rm BG}_k,G^{\rm BG}_k,Q^{\rm BG}_k$ from harmonic forms on a Poincar\'e-Einstein manifold with conformal infinity $(M,[h_0])$.
We say that a $k$-form $\omega$ is \emph{polyhomogeneous} on $\bar{X}$ if it is smooth on $X$ and with an expansion at the boundary $M=\{x=0\}$
\[\omega\sim \sum_{j=0}^{\infty}\sum_{\ell=0}^{\ell(j)}x^j\log(x)^\ell
(\omega_{j,\ell}^{(t)}+\omega_{j,\ell}^{(n)}\wedge dx)\]
for some forms $\omega_{j,\ell}^{(t)}\in C^{\infty}(M,\Lambda^{k}(M))$ and $\omega^{(n)}_{j,\ell}\in C^{\infty}(M,\Lambda^{k-1}(M))$ and some sequence $j\in\nn_0\to \ell(j)\in\nn_0$. 
We show that the Branson-Gover operators appear naturally in the resolution of the absolute or relative Dirichlet type problems for the Laplacian on forms on $\bar{X}$.
\begin{theo}\label{theo2}
Let $(X^{n+1},g)$ be an odd-dimensional Poincar\'e-Einstein manifold with conformal infinity $(M,[h_0])$, let 
$k<\ndemi$ and $\alpha\in(0,1)$.\\
\textbf{(i)} For any $\omega_0\in C^{\infty}(M,\Lambda^{k}(M))$, harmonic
forms $\omega\in C^{\ndemi-k,\alpha}(\bar{X},\Lambda^k(\bar{X}))$ with boundary value $\omega|_{M}=\omega_0$ exist, are unique modulo $\ker_{L^2}(\Delta_k)$ and are actually polyhomogeneous with an expansion at $M$ at order $O(x^{n-2k+1})$ given by 
\[\begin{gathered}
\omega=\omega_0+\sum_{j=1}^{\ndemi-k}x^{2j}(\omega_{j}^{(t)}+\omega_j^{(n)}\wedge \frac{dx}{x})+x^{n-2k}\log(x)L_k\omega_0\\
+x^{n-2k+1}\log(x)(G_k\omega_0)\wedge dx+O(x^{n-2k+1})
\end{gathered}\]
where $L_k,G_k$ are, up to a normalization constant, the Branson-Gover operators in \eqref{Lk}, \eqref{Gk} and $\omega^{(\cdot)}_j$ are forms on $M$.\\
\textbf{(ii)} For any closed form $\omega_0\in C^{\infty}(M,\Lambda^{k-1}(M))$, harmonic forms $\omega$ such that 
$x\omega\in C^{\ndemi-k+1,\alpha}(\bar{X},\Lambda^{k}(\bar{X}))$ and $\omega=x^{-1}(\omega_0\wedge dx)+O(x)$ exist, are unique modulo 
$\ker_{L^2}(\Delta_k)$ and $x\omega$ is polyhomogeneous with expansion at $M$ given by
\[\begin{gathered}
\omega=\omega_0\wedge \frac{dx}{x}+\sum_{j=1}^{\ndemi-k}x^{2j}({\omega'_j}^{(t)}+{\omega'_j}^{(n)}\wedge\frac{dx}{x})+
x^{n-2k+1}\log(x)(Q_{k-1}\omega_0)\wedge dx+O(x^{n-2k+1})
\end{gathered}\]
where $Q_{k-1}$ is, up to a normalization constant, the operator \eqref{Qk} of Branson-Gover and ${\omega'_j}^{(\cdot)}$ are smooth forms on $M$.
\end{theo}

The Dirichlet problem for functions in this geometric setting is studied by 
Graham-Zworski \cite{GRZ} and Joshi-Sa Barreto \cite{JSB}. In a more general setting (but again for functions), 
it was analyzed by Anderson \cite{An} and Sullivan \cite{Su}.\\

We also prove in Subsection \ref{qcurvature} that, with $Q_0$ defined by the Theorem above,
\[Q_01=\frac{n(-1)^{\ndemi+1}}{2^{n-1}\ndemi!(\ndemi-1)!}Q\]
where $Q$ is Branson $Q$-curvature. So $Q$ can be seen as an obstruction to find a harmonic 
$1$-form $\omega$ with $x\omega$ having a high regularity at the boundary and value $dx$ 
at the boundary.\\

In addition, this method allows to obtain the conformal change law of $L_k,G_k,Q_k$, 
the relations between these operators, and some of their analytic properties (e.g. symmetry of $L_k$ and $Q_k$)
 see Subsection \ref{propoflk} and Section \ref{qcurvature}.\\

Next, we analyze the set of regular closed and coclosed forms on $\bar{X}$. Recall that
on a compact manifold $\bar{X}$ with boundary, equipped with a smooth metric $\bar{g}$, there is 
an isomorphism  
\[H^{k}(\bar{X})\simeq \{\omega\in C^{\infty}(\bar{X},\Lambda^k(\bar{X})); d\omega=\delta_{\bar{g}}\omega=0, i_{\pl_n}\omega|_{\pl\bar{X}}=0\}\]
where $\pl_n$ is a unit normal vector field to the boundary, and the absolute cohomology $H^{k}(\bar{X})$ is $\ker d/{\rm Im }\,d$ where $d$ acts on smooth forms. Moreover, one has the long exact sequence in cohomology
\begin{equation}\label{longex}
\dots \xrightarrow{}H^{k-1}(\pl\bar{X})\xrightarrow{} H^k(\bar{X},\pl\bar{X})\xrightarrow{}H^k(\bar{X})\xrightarrow{}H^{k}(\pl\bar{X})\xrightarrow{}H^{k+1}(\bar{X},\pl\bar{X})\xrightarrow{}\dots
\end{equation}
and all these spaces are represented by forms which are closed and coclosed, the maps in the sequence are canonical with respect to $\bar{g}$.
In our Poincar\'e-Einstein case $(X,g)$, say when $k<\ndemi$, only the space $H^k(\bar{X},\pl\bar{X})$ in the long exact sequence 
has a canonical basis of closed and coclosed representatives with respect to $g$ (the $L^2$ harmonic forms), 
in particular there is no canonical metric on the boundary induced by $g$ but  only a canonical conformal class.  
We prove
\begin{theo}\label{th3}
Let $(X^{n+1},g)$ be an odd dimensional Poincar\'e-Einstein manifold with conformal infinity $(M,[h_0])$ and let $k\leq \ndemi$.
Then the spaces 
\[Z^k(\bar{X}):=\{\omega\in C^{n-2k+1}(\bar{X},\Lambda^{k}(\bar{X})); d\omega=\delta_g\omega=0\}\]
are finite dimensional and, if the obstruction tensor of $[h_0]$ vanishes, they are equal to $\{\omega\in C^{\infty}(\bar{X},\Lambda^{k}(\bar{X})); d\omega=\delta_g\omega=0\}$. Then, we have\\ 
\textbf{(i)} For $k<\ndemi$ there is a canonical exact sequence  
\[ 0\xrightarrow{}H^{k}(\bar{X},M)\xrightarrow{} Z^k(\bar{X})\xrightarrow{}
\mc{H}^k(M)\xrightarrow{}H^{k+1}(\bar{X},M)  \]
where $\mc{H}^k(M)$ is the set of conformal harmonics defined in \eqref{hkhkl}.\\
\textbf{(ii)} Let $[Z^k(\bar{X})]$ and $[\mc{H}^k(M)]$ be respectively the image of $Z^k(\bar{X})$ and $\mc{H}^k(M)$ by the natural cohomology maps $Z^k(\bar{X})\to H^k(\bar{X})$
and $\mc{H}^k(M)\to H^k(M)$. Then there is a canonical complex with respect to $g$
\begin{equation}\label{longcomplex}
0\xrightarrow{}\dots \xrightarrow{\iota^k} [Z^k(\bar{X})]\xrightarrow{r^k}
[\mc{H}^k(M)]\xrightarrow{d^k_e}H^{k+1}(\bar{X},M)\xrightarrow{\iota^{k+1}}
[Z^{k+1}(\bar{X})]\xrightarrow{}\dots\xrightarrow{}H^{\ndemi}(\bar{X},M)
\end{equation}
whose cohomology vanishes except possibly the spaces $\ker \iota^k/{\rm Im}\,d^{k-1}_e$.\\ 
\textbf{(iii)} $[\mc{H}^k(M)]= H^k(M)$ if and only if $[Z^k(\bar{X})]=H^k(\bar{X})$ and $\ker\iota^{k+1}={\rm Im}\, d^k_e$. If this holds for all $k\leq \ndemi$
this is a canonical realization of (half of) the long exact sequence \eqref{longex} with respect to $g$.
\end{theo}

The surjectivity of the natural map $\mc{H}^k(M)\to H^{k}(M)$ is named \emph{$(k-1)$-regularity} by Branson and Gover, while
\emph{$(k-1)$-strong regularity} means that the map is an isomorphism, or equivalently $\ker L_{k-1}=\ker d$ (see 
\cite[Th.2.6]{BG}). Thus, $(k-1)$ regularity means that the cohomology group can be represented by conformally invariant representatives.
If $H^{k+1}(X,M)=0$, our result implies that $(k-1)$-regularity means that 
the absolute cohomology group $H^k(\bar{X})$ can be represented by $C^{n-2k+1}(\bar{X},\Lambda^k(\bar{X}))$ forms in $\ker d+\delta_g$. 
We give a criteria for $(k-1)$-regularity:
\begin{prop}\label{qpositif}
Let $(M,[h_0])$ be a compact conformal manifold. If $Q_{k}$ is a positive operator on closed forms in the sense that
$\cjg Q_k\omega,\omega\cjd_{L^2}\geq 0$ for all $\omega\in C^\infty(M,\Lambda^k(M))\cap \ker d$, then $\mc{H}^k(M)\to H^k(M)$ is surjective.
\end{prop}

We should also remark that $(k-1)$-regularity holds for all $k=1,\dots,\ndemi$ if for instance $(M,[h_0])$ contains an Einstein metric in $[h_0]$, this is a result of Gover and Silhan \cite{GS}. If $n=4$, $L_{\ndemi-2}=L_0$ is the Paneitz operator (up to a constant factor) and using a result of Gursky \cite{Gur}, we deduce that if the Yamabe invariant $Y(M,[h_0])$ is positive and  
\[\int_M Q{\rm dvol}_{h_0}+\frac{1}{6}Y(M,[h_0])^2>0\]
then $\mc{H}^1(M)\simeq H^1(M)$ and there is a basis of conformal harmonics of $H^1(M)$.\\

\textbf{Ackowledgement} This work is dedicated to Tom Branson, unfortunetly we could not finish
it before the special volume of SIGMA in his honour appeared. 
We thank Rafe Mazzeo for suggesting to find canonical representatives in $H^k(\bar{X})$ with respect to $g$. 
We also thank Rod Gover for discussions about his paper with Tom Branson.
C.G. is supported by NSF grant DMS0500788, and ANR grants  ANR-05-JCJC-0107091 and 05-JCJCJ-0087-01. 

\section{Poincar\'e-Einstein manifolds and Laplacian on forms} 

\subsection{Poincar\'e-Einstein manifolds}
Let $(X,g)$ be a Poincar\'e-Einstein manifold with conformal infinity $(M,[h])$.
Graham-Lee and Graham \cite{GRL,GR} proved that for any conformal representative $h_0\in[h]$, there exists 
a boundary defining function $x$ of $M=\pl\bar{X}$ in $\bar{X}$ such that 
\[|dx|^2_{x^2g}=1 \textrm{ near }\pl\bar{X}, \quad x^2g|_{TM}=h_0,\]
moreover $x$ is the unique defining function near $M$ satisfying these conditions. 
Such a function is called a \emph{geodesic boundary defining function}
and if $\psi$ is the map $\psi:[0,\eps]\x M\to \bar{X}$ defined by $\psi(t,y):=\psi_t(y)$ where $\psi_t$ is the flow 
of the gradient $\nabla^{x^2g}x$, then $\psi$ pulls the metric $g$ back to 
\[\psi^*g=\frac{dt^2+h_t}{t^2}\] 
for some one-parameter family of metrics on $M$ with $h_0=x^2g|_{TM}$. 
In other words the special form \eqref{modelform} of the metric near infinity 
is not unique and correspond canonically to a geodesic boundary defining function, or equivalently to a conformal 
representative of $[h_0]$.\\

We now discuss the structure of the metric near the boundary, the reader can refer
to Fefferman-Graham \cite[Th 4.8]{FG} for proofs and details.
Let us define the endomorphism $A_x$ on $TM$ corresponding to $\pl_xh_x$ with respect to $h_x$, i.e. as matrices  
\[A_x=h_x^{-1}\pl_x h_x.\]
Then the Einstein condition $\Ric(g)=-ng$ is equivalent to the following differential equations on $A_x$
\[\displaylines{x\pl_xA_x+(1-n+\frac{x}{2}\tra(A_x))A_x=2xh_x^{-1}\Ric(h_x)+\tra(A_x){\rm Id}\cr
\delta_{h_x}(\pl_x h_x)=d\tra(A_x)\cr
\pl_x\tra(A_x)+\frac{1}{2}|A_x|^2=\frac{1}{x}\tra(A_x)}\]
A consequence of these equations and \eqref{hx} is that if $\textrm{Ric}(g)=-ng+O(x^{n-2})$, then 
$h_x$ has an expansion at $x=0$ of the form 
\[h_x=\left\{\begin{array}{ll}
h_0+\sum_{j=1}^{\ndemi-1}x^{2j}h_{2j} +h_{n,1}x^n\log x+O(x^n)& \textrm{if }n\textrm{ is even}\\
h_0+\sum_{j=1}^{(n-1)/2}x^{2j}h_{2j}+O(x^{n}) &\textrm{if }n\textrm{ is odd}
\end{array}\right.
\]
for some tensors $h_{2j}$ and $h_{n,1}$ on $M$, depending in a  natural way on $h_0$ and covariant derivatives of its
Ricci tensor. When $n$ is even, the tensor $h_{n,1}$ is the \emph{obstruction tensor} of $h_0$ in the terminology of 
Fefferman-Graham \cite{FG}, it is trace free (with respect to $h_0$) and so the first log term in 
$A_x$ is $nh_0^{-1}h_{n,1}x^{n-1}\log(x)$. 
A smooth Poincar\'e-Einstein manifold such that $h_x$ has only even powers of $x$ in the Taylor expansion at $x=0$
is called an \emph{smooth even Poincar\'e-Einstein manifold}. 
If $n$ is even and $h_{n,1}=0$, the metric $h_x$ is a smooth even Poincar\'e-Einstein manifold.
When $n$ id odd, the term $\pl_x^nh_x|_{x=0}$ is trace free with respect to $h_0$, which implies 
that $A_x$ has an even Taylor expansion at $x=0$ to order $O(x^{n-1})$. If $\pl^n_xh_x|_{x=0}=0$, then
$h_x$ has an even Taylor expansion in powers of $x$ at $x=0$ with all 
coefficients formally determined by $h_0$. The equations satisfied by $A_x$ easily give (see \cite{FG}) the first terms in the expansion
\begin{equation}\label{firstterms}
h_x=h_0-x^2\frac{P_{0}}{2}+O(x^4), \quad 
\textrm{ where }P_0=\frac{1}{n-2}\Big(2\,\Ric_0-\frac{{\rm Scal}_0}{n-1}h_0\Big),
\end{equation}
$P_0$ is the Schouten tensor of $h_0$, $\Ric_0$ and ${\rm Scal}_0$ are the Ricci and scalar curvature of $h_0$.
 
\subsection{The Laplacian, $d$ and $\delta$}
Let $\Lambda^k(\bar{X})$ be the bundle of $k$-forms on $\bar{X}$. Since for the problem we consider it is somehow quite natural, we will also use along the paper the $b$-bundle of $k$-forms
on $\bar{X}$ in the sense of \cite{Mel}, it will be denoted $\Lambda^k_b(\bar{X})$. This is the exterior product of the $b$ cotangent bundle $T_b^*\bar{X}$, which is canonically isomorphic to $T^*\bar{X}$ over the interior $X$ and whose local basis near a point of the boundary $\pl\bar{X}$ is given by $dy_1,\dots,dy_n,dx/x$ where $y_1,\dots,y_n$ are local coordinates on $\pl\bar{X}$ near this point. We refer the reader to Chapter 2 of \cite{Mel} for a complete analysis about $b$-structures. 
Of course one can pass from $\Lambda^k(\bar{X})$ to $\Lambda^k_b(\bar{X})$ obviously when considering forms on $X$.  
The restriction $\Lambda_b^k(U_\eps)$ of $\Lambda^k(\bar{X})$ to the collar neighbourhood $U_\eps:=[0,\eps]\x M$ of $M$ in $\bar{X}$ can be decomposed as the direct sum 
\[\Lambda_b^k(U_\eps)=\Lambda^k(M)\oplus (\Lambda^{k-1}(M)\wedge \frac{dx}{x})=:\Lambda_t^{k}\oplus \Lambda^k_n.\]
In this splitting, the exterior derivative $d$ and its adjoint $\delta_g$ with respect to $g$ have the form
\begin{equation}\label{formule}
d=\left(
\begin{matrix}
d & 0 \\
(-1)^kx\pl_x & d 
\end{matrix}\right), \quad \delta=\left(
\begin{matrix}
x^2\delta_{x} & (-1)^{k}\star_{x}^{-1}x^{-2k+n+3}\pl_xx^{2k-n-2}\star_x \\
0 & x^2\delta_{x} 
\end{matrix}\right)\end{equation}
and the Hodge Laplace operator is given by  
\begin{equation}\label{decompP}
\begin{gathered}
\Delta_k=\left(
\begin{matrix}
-(x\pl_x)^2+(n-2k)x\pl_x& 2(-1)^{k+1}d \\
0 & -(x\pl_x)^2+(n-2k+2)x\pl_x\end{matrix}
\right)\\
+ \left(\begin{matrix}
x^2\Delta_x -x\star_x^{-1}[\pl_x,\star_x]x\pl_x& (-1)^{k}x\bigl[d,\star_x^{-1}[\pl_x,\star_x]\bigr] \\
2(-1)^{k+1}x^2\delta_x+(-1)^kx^3\bigl[\star^{-1}_x[\pl_x,\star_x],\delta_x\bigr] 
& x^2\Delta_x-x\pl_xx\star^{-1}_x[\pl_x,\star_x]\end{matrix}
\right)\\
=P+P'.
\end{gathered}
\end{equation}
where here, the subscript $\cdot_x$ means ``with respect to the metric $h_x$ on $M$'' and $d$ in the matrices is the exterior derivative on $M$.
Note that $P$ is the indicial operator of $\Delta_k$ in the terminology of \cite{Mel}.

If $H$ is an endomorphism of $TM$, we denote $J(H)$ the operator 
on $\Lambda^k(M)$
\begin{equation}\label{jh}
J(H)(\alpha_1\wedge\dots\wedge\alpha_k):=\sum_{i=1}^k\alpha_1\wedge\dots\wedge\alpha_i(H)\wedge\dots\wedge\alpha_k.
\end{equation}
When $H$ is symmetric, a straightforward computation gives $\star_0J(H)+J(H)\star_0=\tra (H)\star_0$ and so 
\begin{equation}\label{anticom}
[\star_0, J(H)]=2\star_0J(H)-\tra (H)\star_0\end{equation}
Let us define the following operators on $k$-forms on $M$
\begin{equation}\label{defric}
A=J(h_0^{-1}P_0)-\frac{\tra(h_0^{-1}P_0)}{2}{\rm Id}=\frac{2J(h_0^{-1}\Ric)}{n-2}-\frac{n+2k-2}{2(n-1)(n-2)}\Scal_0{\rm Id}.
\end{equation}
Using the approximate Einstein equation for $g$, we obtain
\begin{lem}\label{formlapl}
The operator $\Delta_k$ has a polyhomogeneous expansion at $x=0$ and the first terms in the expansion are given by
\begin{equation}\label{lapgrad}
\begin{gathered}
\Delta_k=P
 +x^2 \left(\begin{matrix}
\Delta_0-x\pl_xA & (-1)^k[d,A]\\
 2(-1)^{k+1}\delta_0& \Delta_0-(2+x\pl_x)A\end{matrix}
\right)
+\sum_{i=2}^{[\ndemi]}x^{2i}\left(\begin{matrix}
R_i+P_ix\pl_x & Q_i\\
Q'_i& R'_i+P'_ix\pl_x\end{matrix}
\right)\\
+nx^n\log(x)
\left(\begin{matrix}
J(h_0^{-1}h_{n,1})x\pl_x & (-1)^{k+1}[d,J(h_0^{-1}h_{n,1})]\\
0&J(h_0^{-1}h_{n,1})(n+x\pl_x)\end{matrix}\right)+O(x^n)
\end{gathered}
\end{equation} 
where $A$ is defined in \eqref{defric} and where the operators $P_i,P'_i,Q_i, Q'_i, R_i$ and $R'_i$ are universal differential operators on $\Lambda(M)$ that can be expressed in terms of covariant derivatives of the Ricci tensor of $h_0$. Moreover the operators $R_i$ and $R'_i$ are of order at most $2$, the $Q_i,Q'_i$ are of order at most $1$ and the $P_i,P'_i$ are of order $0$. If $k=0$, the $x^n\log(x)$ coefficient vanishes.
Finally, if $(X,g)$ is smooth Poincar\'e-Einstein, then $\Delta_k$ is a smooth differential operator on $\bar{X}$, and 
if $(X,g)$ is smooth even Poincar\'e-Einstein, then $\Delta_k$ has an even expansion. 
\end{lem}
\textsl{Proof}: The polyhomogeneity comes from that of the metric $g$. It is moreover 
a smooth expansion if $x^2g$ is smooth on $\bar{X}$. 
A priori, by \eqref{decompP} the first $\log x$ term in the expansion of $\Delta$ at $x=0$ appear at order (at least) $x^n\log x$ and it comes from the diagonal terms in $P_3$ in \eqref{formule}.
Let us define $p=[\ndemi]$ so that the metric $h_x$ has even powers in its expansion at $x=0$ 
up to order $x^{2p+1}$.
We set $D$ the Levi-Civita connexion of the metric $x^2g=dx^2+h_x$. Since $D_{\pl_x}\pl_x=0$ and $D_{\pl x}\pl_{y_i}=\frac{1}{2}\sum_{jk}\pl_xh_{ij}h^{kj}\pl_{y_k}$, the matrix $O_x$ of the parallel transport along the geodesic $x\mapsto (x,y)$ (with respect to the basis $(\pl_{y_i})$) satisfies $D_{\pl_x}O_x(\pl_{y_i})=0$, hence $\pl_x O_x=-\frac{1}{2}A_x\times O_x$ where $A_x$ is the endomorphism $h^{-1}_x\pl_xh_x$. Note that $A_x$ has a Taylor expansion with only odd powers of $x$ up to $x^{2p}$ and the first log term is $nh_0^{-1}h_{n,1}x^{n-1}\log(x)$. 
We infer that $O_x$ is polyhomogeneous in the $x$ variable and has only even powers of $x$ in its Taylor expansion up to $x^{2p}$, the first log term is $-\frac{h_0^{-1}h_{n,1}}{2}x^{n}\log(x)$.
By \eqref{firstterms}, we have $\pl^2_xh|_{x=0}=-P_0$, hence
\[A_x=-xh_0^{-1}P_0+O(x^2), \quad O_x=\textrm{Id}+\frac{1}{4}x^2h_0^{-1}P_0+O(x^3).\]
We note also $O_x$ the parallel transport map.
Now the operator $I_x(\alpha_{1}\wedge\cdots\wedge\alpha_{k})=\alpha_{1}(O_x)\wedge\cdots\wedge \alpha_{k}(O_x)$ is an isometry from $\Lambda^k(M,h_x)$ to $\Lambda^k(M,h_0)$. So we have $\star_x=I^{-1}_x\star_0 I_x$ and we infer that $\star_x$ itself is an operator with a polyhomogeneous expansion in $x$ and with 
only even powers of $x$ in its taylor expansion up to $x^{2p}$, the first log term being 
$\frac{1}{2}x^n\log(x)[J(h_0^{-1}h_{n,1}),\star_0]=-x^n\log(x)\star_0J(h_0^{-1}h_{n,1})$ 
by \eqref{anticom}. 
Since we have 
$$\displaylines{\hfill[\pl_x,\star_x]=\pl_x(\star_x),\hfill\pl_x(\star_x)|_{x=0}= [\star_0,\pl_xI_x|_{x=0}]=0\hfill\mbox{and}\hfill\pl_x^2(\star_x)|_{x=0}=[\star_0,\pl_x^2I_x|_{x=0}]\hfill}$$ we get that $[\pl_x,\star_x]$ is 
polyhomogeneous with  only odd powers of $x$ up to order $x^{2p}$, with first log term $-nx^{n-1}\log(x)\star_0J(h_0^{-1}h_{n,1})$, and that
$$\displaylines{[\pl_x,\star_x]=\pl_x(\star_x)=x\star_0\bigl(J(h_0^{-1}P_0)-\frac{\Scal_0}{2(n-1)}{\rm Id}\bigr)+O(x^2).}$$
Since $\delta_x=(-1)^k\star_x^{-1}d\star_x$, the operators $x\delta_x$ and $x^2[\star^{-1}[\pl_x,\star_x],\delta_x]$ are odd in $x$ up to $O(x^{2p+2})$. By the same way, $x^2[d,\star^{-1}[\pl_x,\star_x]]$ is odd up to order $x^{2p+2}$ and the operators $\star_{x}^{-1}[\pl_x,\star_{x}]x(k-x\pl_x)$, $x^2\Delta_x$ and $ (k-\pl_xx)x\star^{-1}[\star_x,\pl_x]$ are even in $x$ up to $O(x^{2p+1})$. This achieves the proof by gathering all these facts.
\qed\\

\subsection{Indicial equations}\label{indeq}
We give the indicial equations satisfied by $\Delta_k$, which are essential to the construction of formal 
power series solutions of $\Delta_k\omega=0$.\\ 

\textbf{Notation}: If $f$ is a function on $\bar{X}$ and $\omega$ a 
$k$-form defined near the boundary, 
we will say that $\omega$ is a $O_n(f)$ (resp. $O_t(f)$) if its $\Lambda^k_n$ (resp. $\Lambda^k_t$) components are $O(f)$.\\

For $\la\in\cc$, the operator $x^{-\la}\Delta_kx^{\la}$ can be considered near the boundary as a family of 
operators on $\Lambda^{k}_t\oplus\Lambda^k_n$ depending on $(x,\la)$, and for any 
$\omega\in C^{\infty}(U_\eps,\Lambda_t^k\oplus\Lambda^k_n)$ one has 
\begin{equation}\label{indic}
x^{-\la}\Delta_k(x^{\la}\omega)=P_\la\Big(\omega_0^{(t)}+\omega_0^{(n)}\wedge \frac{dx}{x}\Big)+O(x)
\end{equation}
where $P_{\la}:=x^{-\la}Px^\la$, $\omega_0^{(t)}=(i_{x\pl_x}(\omega\wedge\frac{dx}{x}))|_{x=0}$ and
$\omega^{(n)}_0:=(i_{x\pl_x}\omega)|_{x=0}$. The operator $P_{\la}$ is named indicial family and is a 
one-parameter family of operators on $\Lambda^k_n\oplus \Lambda^k_t$ viewed as a bundle over $M$, its expression is  
\[P_\la=\left(\begin{matrix}
-\la^2+(n-2k)\la& 2(-1)^{k+1}d \\
0 & -\la^2+(n-2k+2)\la\end{matrix}
\right)\]
The \emph{indicial roots} of $\Delta_k$ are the $\la\in\cc$ such that $P_\la$ is not invertible on 
the set of smooth sections of $\Lambda^k_t\oplus\Lambda^k_n$ over $M$, i.e. on $C^{\infty}(M,\Lambda^k(M)\oplus\Lambda^{k-1}(M))$. 
In our case, a simple computation shows that these are given by $0,n-2k,0,n-2k+2$. The 
first two roots are roots in the $\Lambda^k_t$ component and the last two are roots in the $\Lambda^k_n$ component. 
In particular, this proves that for $j$ not a root, and $(\omega_0^{(t)},\omega_0^{(n)})\in \Lambda^k(M)\oplus\Lambda^{k-1}(M)$,
there exists a unique pair $(\alpha_0^{(t)},\alpha_0^{(n)})\in \Lambda^{k}(M)\oplus\Lambda^{k-1}(M)$
such that near $M$
\[\Delta_k\Big(x^j\alpha_0^{(t)}+x^j\alpha_0^{(n)}\wedge \frac{dx}{x}\Big)=x^j\Big(\omega_0^{(t)}+\omega_0^{(n)}\wedge \frac{dx}{x}\Big)+
O(x^{j+1})\]
More precisely, and including coefficients with $\log$ terms, we have for $l\in\nn^*$ (resp. $l=0$)
\begin{equation}\label{indic0}
\begin{gathered}
\Delta_k x^j\log^l(x)
\left(\begin{array}{l}
\omega_0^{(t)}\\
\omega_0^{(n)}
\end{array}\right)=x^j\log^l(x)\left(
\begin{array}{l}
j(n-2k-j)\omega_0^{(t)}+2(-1)^{k+1}d\omega_0^{(n)}\\
j(n-2k+2-j)\omega_0^{(n)}
\end{array}\right)\\
\quad\quad\quad\quad\quad\quad\quad\quad+O(x^{j}\log^{l-1}(x))\quad\bigr({\rm resp}.~~+O(x^{j+1})\bigl)
\end{gathered}\end{equation}
if $\omega_0^{(t)},\omega_0^{(n)}\in C^{\infty}(M,\Lambda^k(M)\oplus\Lambda^{k-1}(M))$,
and in the critical cases, for any $l\in\nn_0=\{0\}\cup \nn$
\begin{equation}\label{indic1}
\begin{gathered}
\Delta_k (\log^l(x)\omega_0^{(t)})=  \,l(n-2k)\log^{l-1}(x)\omega_0^{(t)}-l(l-1)\log^{l-2}(x)\omega_0^{(t)}+O(x^2\log x)\\
\Delta_k (x^{n-2k}\log^l(x)\omega_0^{(t)})= 
\, l(2k-n)x^{n-2k}\log^{l-1}(x)\omega_0^{(t)}-l(l-1)x^{n-2k}\log^{l-2}(x)\omega_0^{(t)}\\
+O(x^{n-2k+2}\log^l(x))\\
\Delta_k \bigl(x^{n-2k+2}\log^l(x)\omega^{(n)}\wedge\frac{ dx}{x}\bigr)=
l(2k-2-n)x^{n-2k+2}\log^{l-1}(x)\omega^{(n)}\wedge\frac{ dx}{x}\\
-l(l-1)x^{n-2k+2}\log^{l-2}(x)\omega^{(n)}\wedge\frac{dx}{x}
+O(x^{n-2k+3}\log^l(x)).
\end{gathered}
\end{equation}

\section{Absolute and relative Dirichlet problems} 

The goal of this section is to solve the Dirichlet type problems for $\Delta_k$ when $k<\ndemi$ for the two natural boundary conditions. Note that the vector field $x\pl_x$ can be seen as the unit, normal, inward vector field to $M$ in $\bar{X}$. A $k$-form $\omega\in\Lambda^k_b(\bar{X})$ is said to satisfy the absolute (resp. the relative) boundary condition if
$$\displaylines{\hfill\lim_{x\to0}i_{x\pl_x}\omega=0\hfill
({\rm resp.}\quad\lim_{x\to0}i_{x\pl_x}(\frac{dx}{x}\wedge\omega)=0).\hfill}$$

We denote $C^{p,\alpha}(\bar{X},\Lambda_b^k(\bar{X}))$ the sections of $\Lambda_b^k(\bar{X})$ which are $C^{p,\alpha}$, equivalently $i_{x\pl_x}\omega$ and $i_{x\pl_x}\bigl(\frac{dx}{x}\wedge\omega)$ are $C^{p,\alpha}$ on $\bar{X}$.
\subsection{Absolute boundary condition}

\begin{prop}\label{poisson}
Let $k<\ndemi$, $\alpha\in(0,1)$ and $\omega_0\in C^{\infty}(M,\Lambda^k(M))$.\\
\textbf{(i)} There exists a 
solution $\omega$ to the following absolute Dirichlet problem:
\begin{equation}\label{dirichletpb}
\left\{\begin{array}{l}\omega\in C^{n-2k-1,\alpha}(\bar{X},\Lambda^k(\bar{X})),\\
\Delta_k\omega=0\mbox{ on }X,\\
\omega|_{M}=\omega_0,\dis\lim_{x\to0}i_{x\pl_x}\omega=0.\end{array}
\right.
\end{equation}
Moreover, this solution is unique modulo the $L^2$ kernel of $\Delta_k$.\\ 
\textbf{(ii)} The solution $\omega$ is smooth in $\bar{X}$ when $n$ is odd, while it is polyhomogeneous when $n$ is even  
with an expansion at order $x^{n}$ of the form
\begin{equation}\label{omegaexp}
\begin{split}
\omega=&\sum_{j=0}^{n-1}x^j\omega^{(t)}_{j}+\sum_{j=2}^{n-1}x^j\omega^{(n)}_{j}\wedge \frac{dx}{x}+
\log x\Big(\sum_{j=n-2k}^{n-1}x^{j}\omega^{(t)}_{j,1}+\sum_{j=n-2k+2}^{n}
x^{j}\omega^{(n)}_{j,1}\wedge\frac{dx}{x}\Big)\\
 &+\left\{
\begin{array}{ll}
O_t(x^{n}\log x)+O_n(x^{n+1}\log x) & \textrm{ if }k>0\\
O(x^n) & \textrm{ if }k=0
\end{array}\right. 
\end{split}
\end{equation}
as $x\to 0$, where $\omega_{j}^{(\cdot)},\omega_{j,1}^{(\cdot)}$ are smooth forms on $M$. Moreover, we have
\[\begin{gathered}
\omega_{j}^{(t)}=P_j^{(t)}\omega_0 \textrm{ for }j<n-2k,  \quad \omega_j^{(n)}=P_j^{(n)}\omega_0 \textrm{ for }j<n-2k+2\\
\omega_{n-2k,1}^{(t)}=P_{n-2k,1}^{(t)}\omega_0
\end{gathered}\]
where $P_{j}^{(t)},P_j^{(n)},P_{n-2k,1}^{(t)}$ are universal smooth differential operators on $\Lambda(M)$
depending naturally on covariant derivatives of the curvature tensor of $h_0$.\\ 
\textbf{(iii)} If $n$ is even and $(X,g)$ is a smooth Poincar\'e-Einstein manifold, 
then we have $\omega=\omega_1+x^{n-2k}\log(x)\omega_2$ for some forms 
$\omega_1,\omega_2\in C^{\infty}(\bar{X},\Lambda_b^{k}(\bar{X}))$ with $\omega_2=O(x^\infty)$ if and only if 
$\omega_{n-2k,1}^{(t)}=\omega_{n-2k+2,1}^{(n)}=0$.\\
\textbf{(iv)} $\omega$ satisfies $\delta_g\omega=0$. If in addition $\omega_0$ is closed,
then $d\omega\in\ker_{L^2}(\Delta_{k+1})$.
\end{prop}


\subsubsection{Proof of Proposition \ref{poisson}}

To prove this Proposition, we first need a result of Mazzeo \cite{Ma}: 
\begin{theo}[{\bf Mazzeo}]\label{Mazzeo}
For $k<\ndemi$, the operator $\Delta_k$ is Fredholm and 
there exists a pseudodifferential inverse $E$, bounded on $L^2(X)$, 
such that $\Delta_kE=I-\Pi_0$ where $\Pi_0$ is the projection on the finite dimensional 
space $\ker_{L^2}(\Delta_k)$. This implies an isomorphism between $\ker_{L^2}(\Delta_k)$ and the relative
cohomology $H^{k}(\bar{X},\pl\bar{X})$ of $\bar{X}$.
Moreover any $L^2$ harmonic form $\alpha$ is polyhomogeneous with an expansion near $\pl\bar{X}$ of the form 
\begin{equation}\label{alpha}
\alpha\sim 
x^{n-2k}\sum_{j=0}^\infty\sum_{l=0}^{l(j)}(\alpha^{(t)}_{j,l}x^{j}\log(x)^l+x^{j+2}\log(x)^l\alpha^{(n)}_{j,l}\wedge\frac{ dx}{x}) 
\end{equation}
for some $\alpha^{(t)}_{j,l}\in C^{\infty}(M,\Lambda^k(M)),\alpha^{(n)}_{j,l}\in C^{\infty}(M,\Lambda^{k-1}(M))$ and 
some sequence $l:\nn_0\to \nn_0$.
In addition $E$ maps the space $\{\omega\in C^\infty(\bar{X},\Lambda^k(\bar{X})); \omega=O(x^{\infty})\}$ into polyhomogeneous forms on $\bar{X}$ with a behaviour like \eqref{alpha} near $M$.
\end{theo}
\textbf{Remark}: By using duality through the Hodge star  operator $\star_g$, one obtains trivially a corresponding result
for the case $k>\ndemi+1$. In particular, this gives $\ker_{L^2}(\Delta_k)\simeq H^{k}(\bar{X})$ for $k>\ndemi+1$.\\
 
We can precise the second part of this theorem thanks to the indicial identities obtained by \eqref{decompP}.
\begin{cor}
\label{HarmL2}
Any $L^2$ harmonic $k$-form $\alpha$ on $(X,g)$ is polyhomogeneous and has an expansion at order $x^n\log x$ of the form
\[\alpha=x^{n-2k+2}\Bigl(\sum_{j=0}^{n-1}x^j\alpha_{j}^{(t)}+\sum_{j=0}^{n-1}x^j\alpha_{j}^{(n)}\wedge\frac{dx}{x}+O(x^n\log x)\Bigr)\]
where $\alpha_{j}^{(\cdot)}$ are smooth forms on $M$.
If in addition the metric $(X,g)$ is a smooth Poincar\'e-Einstein manifold, 
then $\alpha\in x^{n-2k+2}C^{\infty}(\bar{X},\Lambda^{k}(\bar{X}))$ 
and $E$ maps
\[E:\{\omega\in C^{\infty}(\bar{X},\Lambda^{k}(\bar{X})); \omega=O(x^{\infty}), \Pi_0\omega=0\} \longrightarrow x^{n-2k}C^{\infty}(\bar{X},\Lambda^{k}(\bar{X})).\]
\end{cor}
\textsl{Proof}: 
Note that if 
\[\alpha\sim x^{n-2k}\sum_{j=0}^\infty\sum_{l=0}^{l(j)}\bigl(\alpha^{(t)}_{j,l}x^{j}\log(x)^l+x^{j+2}\log(x)^l\alpha^{(n)}_{j,l}\wedge\frac{dx}{x}\bigr)\quad\mbox{and }\Delta_k\alpha=O(x^\infty),\] 
then the indicial equations in Subsection \ref{indeq} and Lemma \ref{formlapl} imply that $l(0)=0$ and $l(j)\leq 1$ for all $j=1,\dots,n-1$ 
(and for all $j>0$ if $h_x$ is smooth in $x$).
Moreover since $d\alpha=0$ for any $\alpha\in\ker_{L^2}(\Delta_k)$, we first obtain from (\ref{formule})
that $\alpha^{(t)}_{0,0}=0$ and so, by \eqref{indic1} that $l(j)=0$ for all $j=0,\dots,n-1$ (and for all $j>0$ if $h_x$ is smooth). The mapping property of $E$ is straightforward by the same type of arguments 
and the fact that $\Delta_kE\omega=O(x^{\infty})$ for $\omega=O(x^\infty)$ such that $\Pi_0\omega=0$.   
\qed\\


We will now use the relations \eqref{indic}, \eqref{indic0} and \eqref{indic1} to show that the jet of a solution $\omega$ to the Dirichlet problem in Proposition \ref{poisson} is partly determined.
Let $\omega_0\in C^{\infty}(M,\Lambda^k(M))$.
Using (\ref{indic}) and the form \eqref{lapgrad} of $\Delta$, we can construct a smooth form $\omega_{F_1}$ on $\bar{X}$, solution to the problem
\begin{equation}\label{P'}
\left\{\begin{array}{l}
\Delta_k\omega_{F_1}=O_t(x^{n-2k})+O_n(x^{n-2k+2})\\
\omega_{F_1}|_{M}=\omega_0\end{array}\right.
\end{equation}
it can be taken as a polynomial in $x$ 
\begin{equation}\label{omegaF1}
\omega_{F_1}=\sum_{2j=0}^{n-2k-1}x^{2j}\omega^{(t)}_{2j}+\sum_{2l=2}^{n-2k+1}x^{2l}\omega^{(n)}_{2l}\wedge\frac{dx}{x}
\end{equation} 
and it is the unique solution of \eqref{P'} modulo $O_t(x^{n-2k})+O_n(x^{n-2k+2})$.
Moreover, by \eqref{lapgrad} and parity arguments, we see that when $n$ is odd, the remaining term in \eqref{P'} can be repaced by $O_t(x^{n-2k+1})+O_n(x^{n-2k+3})$ (recall also that $h_x$ is smooth in that case).
By construction, the $\omega_j^{(t)},\omega_n^{(n)}$ are forms on $M$ which can be expressed as a differential operators $P_j^{(t)},P_{j}^{(n)}$ on $M$ acting on $\omega_0$, determined by the expansion of $P$ given in \eqref{lapgrad}, i.e. by $h_0$ and the covariant derivatives of its curvature tensor.

The indicial factor in (\ref{indic}) vanishes if and only if $j=n-2k,l=n-2k+2$ and $n$ is even. Therefore,
if $n$ is odd, we can continue the construction and there is a formal series 
\[\omega_\infty=\sum_{j=0}^\infty x^{j}(\omega^{(t)}_j+\omega_j^{(n)}\wedge dx)\]
such that $\Delta_k\omega_\infty=O(x^\infty)$.
The formal form $\omega_\infty$ can be realized by Borel Lemma, in the sense that there exists a form $\omega'_\infty\in C^{\infty}(\bar{X},\Lambda^k(\bar{X}))$ with the same asymptotic expansion than $\omega_\infty$ at all order and then 
$\Delta_k\omega'_\infty=O(x^\infty)$.

Now for $n$ even, we need to add log terms to continue the parametrix:  
by (\ref{indic1}) one can modify $\omega_{F_1}$ to 
\begin{equation}\label{omegaF2}
\omega_{F_2}=\omega_{F_1}+ x^{n-2k}\log(x)\omega^{(t)}_{n-2k,1}
\end{equation} 
such that $\Delta_k\omega_{F_2}= O(x^{n-2k+2}\log x)$. Actually, using \eqref{lapgrad} and parity arguments once more,
we see that
\begin{equation}\label{deltakome}
\Delta_k\omega_{F_2}=2(-1)^{k+1}x^{n-2k+2}\log x(\delta_0\omega^{(t)}_{n-2k,1})\wedge \frac{dx}{x}+O_t(x^{n-2k+2}\log x)+O_n(x^{n-2k+2}).
\end{equation}
Now we want to show
\begin{lem}\label{technic}
The $k$-form $\omega^{(t)}_{n-2k,1}$ on $M$ satisfies $\delta_0\omega_{n-2k,1}^{(t)}=0$.
\end{lem}
\textsl{Proof}: 
From \eqref{deltakome}, and the expression of $\delta$, we obtain \[\delta_g\Delta_k\omega_{F_2}=-2x^{n-2k+2}\,\delta_0\omega_{n-2k,1}+O(x^{n-2k+3}\log x).\]
But $\delta_g\Delta_k\omega_{F_2}=\Delta_{k-1}\delta_g \omega_{F_2}$ and
\[\delta_g\omega_{F_2}=\sum_{j=2}^{n-2k+2}x^{j}\omega'^{(t)}_j+\sum_{l=3}^{n-2k+3}x^l\omega'^{(n)}_l\wedge \frac{dx}{x}
+ x^{n-2k+2}\log(x)\,\delta_0\omega^{(t)}_{n-2k,1}+O(x^{n-2k+3}\log x)\]
for some forms $\omega'^{(.)}_j$ on $M$, so by uniqueness of \eqref{P'} and the fact that $\delta_g\omega_{F_2}=O(x^2)$
we deduce that 
\[\delta_g\omega_{F_2}=x^{n-2k+2}\omega'^{(t)}_{n-2k+2}+ x^{n-2k+2}\log(x)\,\delta_0\omega^{(t)}_{n-2k,1}+O(x^{n-2k+3}\log x).\]
Using now \eqref{indic} and \eqref{indic1}, we obtain $\Delta_{k-1}\delta_g\omega_{F_2}=(2k-n-2)x^{n-2k+2}\delta_0\omega_{n-2k,1}^{(t)}+O(x^{n-2k+3}\log x)$, and since $k<\ndemi$
this implies $\delta_0\omega_{n-2k,1}^{(t)}=0$.
\qed\\

We infer that there is no term of order $x^{n-2k+2}\log x$ in the $\Lambda^k_n$ part of $\Delta_k\omega_{F_2}$ and we can continue to solve the problem modulo $O(x^{\infty})$ using formal power series with log terms using the indicial equations.
The formal solution when $n$ is even will be given by
\begin{equation}\label{formal}
\begin{gathered}
\omega_\infty=\sum_{j=0}^{\frac{n}{2}-1}x^{2j}\omega^{(t)}_{2j}+\sum_{j=1}^{\frac{n}{2}}x^{2j}\omega^{(n)}_{2j}\wedge\frac{dx}{x}
+ \sum_{j=\frac{n}{2}-k}^{\frac{n}{2}-1}x^{2j}\log(x)\omega^{(t)}_{2j,1}\\
+\sum_{j=\frac{n}{2}-k+1}^{\frac{n}{2}}x^{2j}\log(x)\omega^{(n)}_{2j,1}\wedge \frac{dx}{x}+ x^{n}\sum_{j=0}^\infty\sum_{l=0}^{j+1} (\omega_{n+j,l}^{(t)}+x\omega_{n+j,l}^{(n)}\wedge\frac{dx}{x})x^{j}(\log x)^l
\end{gathered}
\end{equation}  
which again is realized through Borel's Lemma to have $\Delta_k\omega_\infty=O(x^\infty)$. 
Notice that when the metric $h_x$ is smooth, the second line 
in \eqref{formal} has $\omega^{(t)}_{j,l}=\omega^{(n)}_{j,l}=0$ for $l>1$ 
since these terms come from the log terms of the expansion of $h_x$ in \eqref{hx} (and thus of $\Delta_{k}$).
The terms $(\omega^{(t)}_{j})_{j<n-2k}$, $(\omega^{(n)}_j)_{j<n-2k+2}$ and $\omega_{n-2k,1}^{(t)}$ are formally determined by $\omega_0$ and are expressed as a differential operator on $M$ acting on $\omega_0$, 
the terms $\omega^{(t)}_{n-2k},\omega^{(n)}_{n-2k+2}$ are formally undetermined, 
the remaining terms are formally determined by $\omega_0,\omega^{(t)}_{n-2k}$ and $\omega^{(n)}_{n-2k+2}$.

So we have proved 
\begin{prop}[{\bf Formal solution}]\label{formalsol}
Let $\omega_0,v^{(t)}, v^{(n)}\in C^\infty(M,\Lambda(M))$, then there exists a form $\omega_\infty\in C^{n-2k-1}(\bar{X},\Lambda_b^k(\bar{X}))$, unique modulo $O(x^\infty)$, which is smooth on $\bar{X}$ when $n$ is odd and with a polyhomogeneous expansion at $\pl\bar{X}$ 
of the form \eqref{formal} when $n$ is even, such that $\Delta_k\omega_\infty=O(x^\infty)$, 
$\omega_\infty|_{\pl\bar{X}}=\omega_0$, $\omega_{n-2k}^{(t)}=v^{(t)}$ 
and $\omega_{n-2k+2}^{(n)}=v^{(n)}$ in the expansion \eqref{formal}. 
\end{prop}

To correct the approximate solution and obtain a true harmonic form, we add $-E\Delta_k(\omega_\infty)$ to $\omega_\infty$ and so
\[\Delta_k(\omega_\infty-E\Delta_k\omega_\infty)=\Pi_0\Delta_k\omega_\infty.\]
We want to prove that $\Pi_0\Delta_k\omega_\infty=0$ or equivalently that $\cjg \Delta_k\omega_\infty,\alpha\cjd=0$ for any 
$\alpha\in\ker_{L^2}(\Delta_k)$.
For that, we use Green's formula on $\{x\geq \varepsilon\}$ and let $\varepsilon\to 0$, together with the asymptotic $\alpha=O(x^{n-2k+1})$ obtained from Theorem \ref{Mazzeo}, $d\alpha=0$ and $\delta\alpha=0$:
\[\int_{x\geq\varepsilon}\langle\Delta\omega_\infty,\alpha\rangle{\rm dvol}_g=(-1)^n\int_{x=\varepsilon}(\star_gd\omega_\infty)\wedge\alpha-(\star_g\alpha)\wedge\delta\omega_\infty
=O(\varepsilon)\to_{\varepsilon\to 0} 0.\]

In view of the mapping properties of $E$ from Theorem \ref{Mazzeo}, we have thus proved that $\omega=\omega_\infty-E\Delta_k\omega_\infty$ is a harmonic $k$-form of $X$ such that $\omega_{|M}=\omega_0$, 
with an asymptotic of the form \eqref{formal} when $n$ is even and smooth on $\bar{X}$ when $n$ is odd,
such that 
\[\omega-\omega_{F_2}=O_t(x^{n-2k})+O_n(x^{n-2k+2})\]
and with $C^{n-2k-1,\alpha}(\bar{X},\Lambda^k(\bar{X}))$ regularity.\\

Let us now consider the problem of uniqueness. If one assumes polyhomogeneity of the solution of 
$\Delta_k\omega=0$ with boundary condition $\omega=\omega_0+o(x)$, the construction above with formal series arguments and indicial equations shows that $\omega$ is unique up to $O(x^{n-2k})$, i.e. the first positive indicial root, 
then of course two such solutions would differ by an $L^2$ harmonic form if $k<\ndemi$. This gives
\begin{lem}\label{uniquephg}
Polyhomogeneous forms satisfying $\Delta_k\omega$ and $\omega=\omega_0+o(x)$ are unique modulo the $L^2$ kernel 
of $\Delta_k$.
\end{lem}

Here, since we want a sharp condition on regularity for uniqueness, i.e. we do not assume polyhomogeneity but $C^{n-2k-1,\alpha}$ regularity, we first need a 
preliminary result.  Let $H^{s}(\Lambda^{k}(M))$ be the Sobolev space of order $s\in\zz$ with $k$-forms values, which we will also denote by $H^{s}(M)$ to simplify. The sections of the bundle $\Lambda^k_t\oplus\Lambda^k_n$ over $M$  
are equipped with the natural Sobolev norm $||.||_{H^{s}(M)}$ induced by $H^s(\Lambda^k(M)\oplus\Lambda^{k-1}(M))$.
Then it is proved by Mazzeo \cite[Th. 7.3]{Maz} the following property\footnote{Notice that the result of Mazzeo
is stated for $0$-elliptic operators with smooth coefficients and acting on functions, 
but it is straightforward to check that it applies on bundles and with polyhomogeneous coefficients like this is the case for even-dimensionnal Poincar\'e-Einstein manifolds.} 
\begin{lem}[\textbf{Mazzeo}]\label{maz2}
Let $k<\ndemi$ and let $\omega\in x^{\alpha}L^2(\Lambda^k(X),\textrm{dvol}_g)$ with $\alpha<-\ndemi$ such that 
$\Delta_k\omega=0$, then for all $N\in\nn$, there exist some 
forms $\omega_{j,l}^{(t)},\omega_{j,l}^{(n)}\in H^{-N}(M)$ for $j,l\in\nn_0$ and some sequence $l:\nn_0\to \nn_0$ 
such that 
\begin{equation}\label{weakexp}
\Big|\Big|\,\omega-\sum_{j=0}^{N-3}\sum_{l=0}^{l(j)}x^j(\log x)^l(\omega_{j,l}^{(t)}-
\omega_{j,l}^{(n)}\wedge \frac{dx}{x})\Big|\Big|_{H^{-N}(M)}=O(x^{N-2-\varepsilon})
\end{equation}
for all $\varepsilon>0$.
\end{lem}
Let $\omega,\omega'$ be two harmonic forms which are $C^{n-2k-1,\alpha}(\bar{X},\Lambda^{k}(\bar{X}))$ and which coincide on the boundary, we want to show that their Taylor expansion at $x=0$ coincide to order $n-2k-1$.
Using Lemma \ref{maz2} with $N$ large enough, we see that the arguments used above on formal series 
(based on the indicial equations) also apply by considering norms $||\cdot||_{H^{-N}(M)}$ on $\Lambda^k_t\oplus\Lambda^k_n$, 
in particular that $l(j)=0$ for $j=0,\dots,n-2k-1$ in \eqref{weakexp} for both $\omega$ and $\omega'$, and that their 
coefficients of $x^j$  for $j=0,\dots,n-2k-1$ in the weak expansion \eqref{weakexp} are the same
for $\omega$ and $\omega'$, these are given by
$\omega^{t}_{j,0}=P_j^{(t)}\omega_0$ and $\omega_j^{(n)}=P^{(n)}_j\omega_0$ (and are then continuous on $M$ since
$\omega\in C^{n-2k+1}(\bar{X},\Lambda^k(\bar{X}))$).
But by uniqueness of the expansion \eqref{weakexp} and the regularity assumption on $\omega,\omega'$, 
this implies that $\omega_j^{(t)},\omega^{(n)}_j$ are the coefficients in the Taylor expansion of 
both $\omega$ and $\omega'$ to order $n-2k-1$. The extra H\"older regularity then gives
that $||\omega-\omega'||_{L^\infty(M)}=O(x^{n-2k-1+\alpha})$, but then this implies that 
$\omega-\omega'\in\ker_{L^2}(\Delta_k)$ thus it is
in the $L^2$ kernel of $\Delta_k$, so our construction is unique modulo $\ker_{L^2}(\Delta_k)$. 
This ends the proof of the solution of \eqref{dirichletpb}.\\

Now to deal with (iv), we notice that $d\omega$ is solution of the problem 
\eqref{dirichletpb} for $(k+1)$-forms with the additional condition that the boundary value is 
$d\omega_0=0$. Note that this requires a priori that $k+1\not=\ndemi$. However, the discussion below in 
Subsection \ref{n/2} about the solutions of $\Delta_{\ndemi}\omega=O(x^\infty)$ gives the same result, namely 
that $d\omega\in\ker_{L^2}(\Delta_{\ndemi})$ if $\omega$ is a solution of \eqref{dirichletpb} with $k=\ndemi-1$.
\qed\\

We conclude this section by a remark.
\begin{prop}\label{propr}
The forms $\omega_{F_1}$ of \eqref{P'} and $\omega$ of Proposition \ref{poisson} satisfy
\[\delta_g\omega=0, \quad \delta_g \omega_{F_1}=O_t(x^{n-2k+2})+O_n(x^{n-2k+4}).\]
\end{prop}
\textsl{Proof}:
Let $\omega$ be the exact solution of $\Delta_k\omega=0$,  $\omega|_{x=0}=\omega_0$ in Proposition \ref{poisson}. Since $\delta_g\Delta_k=\Delta_{k-1}\delta_g$, we deduce that $\omega':=\delta_g\omega$ is solution of $\Delta_k\omega'=0$
with $\omega'|_{x=0}=0$ and moreover it is polyhomogeneous since $\omega$ is polyhomogeneous, 
so Proposition \ref{poisson} and Lemma \ref{maz2} imply that $\delta_g\omega\in\ker_{L^2}(\Delta_{k-1})$ and thus $\delta_g\omega =O(x^{n-2k+3})$ by Corollary \ref{HarmL2}. 
Hence $\delta_g\omega$ is closed and integration by parts on $\{x\geq \eps\}$ shows, by letting $\eps\to 0$, that $\cjg\delta_g\omega,\delta_g\omega\cjd=0$. The part with $\omega_{F_1}$ is also based on $\delta_g\Delta_k=\Delta_{k-1}\delta_g$ and the uniqueness of the solution of \eqref{P'} up to $O_t(x^{n-2k+2})+O_n(x^{n-2k+4})$ on $(k-1)$-forms.
\qed

\subsubsection{The case $k=\frac{n}{2}$}\label{n/2}
 
In this case one only intend to solve the equation $\Delta_k\omega=O(x^{\infty})$, say in the
set of almost bounded forms ($\log x$ times bounded).
The indicial equation tells us that $0$ is a double indicial root for  
the $\Lambda^{k}_t$ part, while $0,2$ are the two simple roots for the $\Lambda_n$ part.
So for $\omega_0,\omega_1\in \Lambda^{k}(M)$, one can construct a polyhomogeneous form 
\begin{equation*}
\begin{split}
\omega_F=\,&\omega_1 \log(x) + \omega_0+\frac{(-1)^{\frac{n}{2}+1}}{2}x^2(\log x)^2\,\delta_0\omega_1\wedge\frac{ dx}{x}\\
&+(-1)^\frac{n}{2} x^2\log(x)\,(-\delta_0\omega_0+\frac{1}{2}\delta_0\omega_1)\wedge\frac{dx}{x}+O_t(x^{2}(\log x)^2)+O_n(x^{3}(\log x)^2)
\end{split}
\end{equation*}
such that $\Delta_k\omega_F=O(x^{\infty})$ and it is unique modulo $O(x^\infty)$ if 
the order $x$ coefficient in the $\Lambda_t$ component is assumed to be $0$.
We also recall a result proved by Yeganefar \cite[Corollary 3.10]{Ye}.
\begin{prop}\label{nader}
For an odd dimensional Poincar\'e-Einstein manifold $(X^{n+1},g)$, there is an isomorphism between $\ker_{L^2}(\Delta_{\ndemi})$ and $H^{\ndemi}(\bar{X},\pl\bar{X})$ and between $\ker_{L^2}(\Delta_{\ndemi+1})$ and $H^{\ndemi+1}(\bar{X})$. 
\end{prop}


\subsection{Relative boundary condition}

\begin{prop}\label{construQ} 
Let $0<k\leq\ndemi-1$, $x$ be a geodesic boundary defining function 
and $\omega_0\in C^\infty(M,\Lambda^{k-1}(M))$ be a closed form. Then there exists a unique, modulo 
$\ker_{L^2}(\Delta_k)$, form $\omega$ such that 
\begin{equation}\label{dx/x}
\left\{\begin{array}{l}
\omega\in C^{n-2k}(\bar{X},\Lambda_b^{k}(\bar{X})),\\
\Delta_k \omega=0\mbox{ on }X,\\
\omega|_{M}=0,\dis\lim_{x\to 0}i_{x\pl_x}\omega=\omega_0.
\end{array}\right.\end{equation}
Moreover $\omega$ is closed, smooth on $\bar{X}$ when $n$ is odd, while it is polyhomogeneous when $n$ is even with an expansion at order $O(x^{n-1}\log x)$ of the form
\begin{equation}\label{omdx/x}
\begin{split}
\omega=&\Big(\sum_{j=0}^{n-1}x^j\omega_{j}^{(n)}\wedge\frac{ dx}{x}+\sum_{j=1}^{n-2}x^j\omega_{j}^{(t)}
\Big)\\
&+x^{n-2k+2}\log x\Big(\sum_{j=0}^{2k-3}x^j\omega_{j,1}^{(n)}\wedge \frac{dx}{x}+\sum_{j=0}^{2k-4}x^j\omega_{j,1}^{(t)}\Big)+O(x^{n-1}\log x)
\end{split}
\end{equation} 
for some forms $\omega_{j}^{(.)},\omega_{j,1}^{(.)}$ on $M$.
\end{prop}
\textsl{Proof}: the proof is similar to that of Proposition \ref{poisson}, so we do not give the full details but we shall use the same notations. We search a formal solution $\omega'_{\infty}$ of $\Delta_k\omega'_{\infty}=0$ with $\omega'_{\infty}=\omega_0\wedge\frac{dx}{x}+O(x)$.
Using the indicial equations in Subsection \ref{indeq} and the form of $\Delta_k$ in Lemma \ref{formlapl}, we can construct the exponents in the formal series as long as the exponent is not a solution of the indicial equation. Since $d\omega_0=0$ by assumption, we have
\[\Delta_k(\omega_0\wedge \frac{dx}{x})=2(-1)^{k+1}d\omega_0+O(x^2)=O(x^2)\]
and so we can continue the construction of $\omega'_{\infty}$ until the power $x^{n-2k}$ in the tangential part $\Lambda_t^k$ and $x^{n-2k+2}$ in the $\Lambda_n^k$ part. At that point, since $x^{n-2k}$ and $x^{n-2k+2}$ are solution of the indical equation of $\Delta_k$ in respectively the $\Lambda^k_t$ and $\Lambda_n^k$ part, there is a 
$x^{n-2k}\log x$ term to include in the $\Lambda_t^k$ part. Using in addition that $\Delta_k$ begins with a sum of even powers of $x$, we see like in Proposition \ref{poisson} that when $n$ is odd, a formal series $\omega'_{\infty}$ with no log terms can be constructed
to solve $\Delta_k\omega'_{\infty}=O(x^{\infty})$, while when $n$ is even we can first construct 
\begin{equation}\label{F2F1}
\omega'_{F_2}=\underbrace{\sum_{2j=0}^{n-2k}x^{2j}\omega_{2j}^{(n)}\wedge\frac{ dx}{x}+\sum_{2j=2}^{n-2k-2}x^{2j}\omega_{2j}^{(t)}}_{=\omega'_{F_1}}+x^{n-2k}\log(x)\,\omega_{n-2k,1}^{(t)}
\end{equation}
with $\omega_0^{(n)}=\omega_0$ so that $\Delta_{k}\omega'_{F_2}=O(x^{n-2k+2}\log x)$, and the coefficients are uniquely determined by $\omega_0$. First observe that $d\omega'_{F_1}=O(x^2)$ satisfies  
$\Delta_{k+1}d\omega'_{F_1}=O_t(x^{n-2k})+O_t(x^{n-2k+2})$ and since the indicial root in $[2,n-2k]$
for $\Delta_{k+1}$ are $n-2k-2$ in the $\Lambda^{k+1}_t$ part and $n-2k$ in the $\Lambda^{k+1}_n$ part, we deduce
that $d\omega'_{F_1}=O_t(x^{n-2k-2})+O_n(x^{n-2k})$ and so 
\begin{equation}\begin{split}\label{domegaF1}
d\omega'_{F_1}=&\,\sum_{2j=2}^{n-2k-2} d\omega^{(t)}_{2j}x^{2j}+
\sum_{2j=2}^{n-2k-2}x^{2j}((-1)^k2j\omega^{(t)}_{2j}+d\omega^{(n)}_{2j})\wedge\frac{dx}{x}+x^{n-2k}d\omega_{n-2k}^{(n)}\wedge \frac{dx}{x}\\
=&\, x^{n-2k}d\omega_{n-2k}^{(n)}\wedge \frac{dx}{x}.
\end{split}
\end{equation}
Now we want to show that $\omega^{(t)}_{n-2k,1}=0$ to continue the construction of the formal solution to higher order.
Clearly now we have 
$$\displaylines{d\omega'_{F_2}=x^{n-2k}\log(x)\bigl(d\omega^{(t)}_{n-2k,1}+(-1)^k(n-2k)\omega^{(t)}_{n-2k,1}\wedge\frac{dx}{x}\bigr)\cr
+x^{n-2k}\bigl(d\omega_{n-2k}^{(n)}+(-1)^k\omega^{(t)}_{n-2k,1}\bigr)\wedge\frac{dx}{x}}$$
and so that
$$\Delta_{k+1} d\omega'_{F_2}=(n-2k)x^{n-2k}\bigl((-1)^{k+1}(n-2k)\omega^{(t)}_{n-2k,1}\wedge\frac{dx}{x}+\log(x)d\omega^{(t)}_{n-2k,1}\bigr)+O(x^{n-2k+1}).$$
But since $d\Delta_k\omega'_{F_2}=O(x^{n-2k+2}\log(x))$, we infer that $\omega_{n-2k,1}^{(t)}$ must vanish,
and we obtain 
\[\Delta_{k}\omega'_{F_2}=\Delta_k\omega'_{F_1}=O(x^{n-2k+2})\]
Since the order $x^{n-2k+2}$ is a solution of the indicial equation in the normal part $\Lambda_n^k$, we need to add a $x^{n-2k+2}\log(x)$ normal term to continue the
construction of the formal solution. Since all the subsequent orders are not solution of the indicial equation for $\Delta_k$, we can construct, using Borel lemma, a polyhomogeneous $k$-form on $X$ with expansion to order $x^{n-1}\log(x)$ of the form given by  \eqref{omdx/x}, which coincides with $\omega_{F_2}$ at order $O_n(x^{n-2k+2}\log x)+O_t(x^{n-2k})$.
To obtain an exact solution of \eqref{dx/x}, we can correct $\omega'_{\infty}$ by setting 
$\omega=\omega'_{\infty}-E\Delta_k\omega'_{\infty}$ where $E$ is defined in Proposition \ref{Mazzeo}.

The argument for the uniqueness modulo $\ker_{L^2}\Delta_g$ is similar 
to that used in the proof of Proposition \ref{poisson}.

To prove that $\omega$ is closed, it suffices to observe that $d\omega\in C^{n-2k-3}(\bar{X},\Lambda_b^{k+1}(\bar{X}))$ and  $d\omega=O(x^2)$ and then use Proposition \ref{poisson} to deduce that $d\omega\in\ker_{L^2}(\Delta_{k+1})$.
Then $\delta_gd\omega=0$ and, considering the decay of $d\omega$ and $\omega$ at the boundary, 
we see by integration by part that $d\omega=0$. 
\qed\\

\noindent\textbf{Remarks}: it is important to remark that the solution $\omega$ of the problem \eqref{dx/x} depends 
on $\omega_0$ but also on the choice of $x$. Note also that the form $\omega$ solution of \eqref{dx/x} satisfies 
$x\omega\in C^{n-2k+1,\alpha}(\bar{X},\Lambda^{k}(\bar{X}))$ for all $\alpha\in(0,1)$.

\section{$L_k$, $G_k$ and $Q_k$ operators}
In this section we suppose that $M$ has an even dimension $n$.
\subsection{Definitions}
The operators $L_k,G_k$ derive from the solution of the absolute Dirichlet problem:
\begin{defi}\label{LkHk}  
For $k<\ndemi$, the operators $L_k: C^\infty(M,\Lambda^k(M))\to C^\infty(M,\Lambda^k(M))$ and $G_k :C^\infty(M,\Lambda^k(M))\to C^{\infty}(M,\Lambda^{k-1}(M))$ are defined by
$L_k\omega_0=:\omega_{n-2k,1}^{(t)}$ and $G_k\omega_0:=\omega_{n-2k+2,1}^{(n)}$
where $\omega_{n-2k,1}^{(t)}, \omega_{n-2k+2,1}^{(n)}$ are given in the expansion \eqref{omegaexp}.
 When $k=\frac{n}{2}$, we define $G_{\frac{n}{2}}:=(-1)^{\ndemi+1}\delta_0$. 
\end{defi}

The operator $Q_k$ derives from the solution of the relative Dirichlet problem:
\begin{defi}\label{definQ}
Let $n$ be even and $k<\ndemi$, the operator $Q_{k-1}:(C^\infty(M,\Lambda^{k-1}(M))\cap \ker d)\to C^{\infty}(M,\Lambda^{k-1}(M))$
is defined by $Q_{k-1}\omega_0:=\omega_{n-2k+2,1}^{(n)}$ where $\omega_{n-2k+2,1}^{(n)}$ is given in the expansion
\eqref{omdx/x}.
\end{defi}

By Corollary \ref{HarmL2}, $L_k$, $G_k$ and $Q_{k}$ do not depend on the choice of the solution $\omega$ in Propositions \ref{poisson} or \ref{construQ}, if $L_k$ depend only on the boundary $(M,[h_0])$,  the operators $G_k$ and 
$Q_k$ may well depend on the whole manifold $(X,g)$ and not only on the conformal boundary. 
We will see that they actually depend only on $(M,[h_0])$ and that they are differential operators.\\

\subsection{A formal construction}

We show that the definition of $L_k,G_k,Q_k$ can be done using only the formal series solutions. Let us first define
\begin{defi}\label{Bk}
For $k<\ndemi$, the operators $B_k,C_k:C^\infty(M,\Lambda^k(M))\to C^{\infty}(M,\Lambda^{k-1}(M))$  
and $D_k:C^\infty(M,\Lambda^{k}(M))\cap\ker d\to C^{\infty}(M,\Lambda^{k}(M))$ are defined by
\begin{equation}
\begin{gathered}\label{defCk}
B_k\omega_0:=\Big(x^{-n+2k-2}i_{x\pl_x}\Delta_k\omega_{F_1}\Big)|_{x=0},\\
C_k\omega_0:=\Big(x^{-n+2k-2}i_{x\pl_x}(\frac{dx}{x}\wedge \delta_g\omega_{F_1})\Big)|_{x=0}\\
D_k\omega_0:=\Big(x^{-n+2k}i_{x\pl_x}d\omega_{F_1})\Big)|_{x=0}
\end{gathered}\end{equation}
where $\omega_{F_1}$ solves \eqref{P'}. 
\end{defi}

\noindent{\bf Remark}: from the indicial equations and Lemma \ref{technic}, $B_k\omega_0$ is $(-1)^k(n-2k+2)$ times the $x^{n-2k+2}\log(x)$ coefficient in 
the $\Lambda^k_n$ part of $\omega'_{\infty}$ defined in Proposition \ref{formalsol} when $v^{(t)}=0$, 
this is a differential operator on $M$ of order $n-2k+1$ since by construction, $\omega_{F_1}$ contains only derivatives of order at most $n-2k-1$ with respect to $\omega_0$. The operator $C_k$
is well defined thanks to Proposition \ref{propr}, and it is a differential operator of order $n-2k$.
As they come from the expansion of $\Delta_k,\delta_g$, they are natural differential operators 
depending only on $h_0$ and the covariant derivatives of its curvature tensor.\\

\subsubsection{The case of $L_k$} It is clear from the proof of Proposition \ref{poisson} that $L_k\omega_0$ is also the coefficient of the $x^{n-2k}\log x$ term in the expansion of $\omega_{F_2}$ defined in \eqref{omegaF2} and of the formal solution $\omega_\infty$ defined in Proposition \ref{formal}.
The indicial equation shows that
\begin{equation}\label{LkF1}
  L_k\omega_0:=\frac{1}{n-2k}\Big(x^{2k-n}i_{x\pl_x}(\frac{dx}{x}\wedge \Delta_k\omega_{F_1})\Big)|_{x=0}
\end{equation}
where $\omega_{F_1}$ solves \eqref{P'}.  

\subsubsection{The case of $G_k$} Let us return to the construction of the formal series solution in the proof of 
Proposition \ref{dirichletpb}. 
Now let $\omega_{F_2}$ defined in \eqref{omegaF2} and 
\[\omega_{F_2}:=\omega_{F_1}+x^{n-2k}v^{(t)}+x^{n-2k}\log(x)\omega_{n-2k,1}^{(t)}\] 
where $v^{(t)}\in C^\infty(M,\Lambda^k(M))$ is an arbitrary form. By construction of $\omega_{F_1},\omega_{F_2}$, 
the fact that $n-2k$ is an indicial root in the $\Lambda^k_t$ component and Lemma \ref{formlapl}, we have 
\[\Delta_k\omega_{F_2}= (-1)^{k+1}x^{n-2k+2}(B_k\omega_0+2\delta_0v^{(t)})\wedge\frac{dx}{x}
+O_t(x^{n-2k+2}\log x)+O_n(x^{n-2k+4}\log x)\]
to solve away the $x^{n-2k-2}$ term in $\Lambda_n^k$ we need to define
\begin{equation}\label{omegaF3}
\omega_{F_3}:=\omega_{F_2}+\frac{(-1)^{k+1}}{n+2-2k}x^{n-2k}\log(x)(B_k\omega_0+2\delta_0v^{(t)})\wedge\frac{dx}{x}
\end{equation}
so that $\Delta_{k}\omega_{F_3}=O_n(x^{n-2k+4}\log(x))+O_t(x^{n-2k+2}\log(x))$. Since $v^{(t)}$ can be chosen arbitrarily,
the coefficient of $x^{n-2k+3}\log(x)$ in the $\Lambda^k_n$ component of the formal solution $\omega_{F_3}$ does not determine a natural operator in term of the initial data $\omega_0$, contrary to the $x^{n-2k}\log(x)$ coefficient in $\Lambda^k_t$. 
In the definition of $G_k$ above, we used an exact solution on $X$ to fix the $v^{(t)}$ term through the Green function, which a priori makes $G_k$ depend on $(X,g)$ and not only on $(M,[h_0])$.
However there is an equivalent way of fixing $\delta_0v^{(t)}$ without solving a global Dirichlet problem
but by adding an additional condition: 
\begin{prop}\label{formalequations}
Let $\omega_0\in C^\infty(M,\Lambda^{k}(M))$, then there is a polyhomogeneous $k$-form $\omega_F$
such that
\begin{equation}\label{condomegaF}
\left\{\begin{array}{l}
\Delta_k\omega_{F}=O_t(x^{n-2k+1})+O_n(x^{n-2k+3})\\
\delta_g\omega_{F}=O(x^{n-2k+3})\\
\omega=\omega_0+O(x)\end{array}\right..
\end{equation}
It is unique modulo $O_t(x^{n-2k})+O_n(x^{n-2k+2})$ and has an expansion of the form
\begin{equation}
\label{omegaF}\begin{split}
\omega_F=\,&\sum_{j=0}^{\frac{n}{2}-k-1}x^{2j}\omega^{(t)}_{2j}+
\sum_{j=1}^{\frac{n}{2}-k}x^{2j}\omega^{(n)}_{2j}\wedge\frac{dx}{x}\\
&+x^{n-2k}\log(x)\Big(L_k\omega_0+x^2\frac{(-1)^{k+1}}{n-2k}(B_k\omega_0-2C_k\omega_0)\wedge\frac{dx}{x}\Big).
\end{split}
\end{equation}
\end{prop}
\textsl{Proof}: First consider the uniqueness. By the discussion above, the condition on $\Delta_k\omega_F$ implies that $\omega_F$ is necessary of the form $\omega_F=\omega_{F_3}$ defined in \eqref{omegaF3} for some $v^{(t)}$.  
Now we notice that $\delta_g\omega_{F_3}=O(x^2)$ satisfies
in particular $\Delta_{k-1}\delta_g\omega_{F_3}=\delta_g\Delta_k\omega_{F_3}=O(x^{n-2k+3})$, and again by the indicial 
equation this implies that $\delta_g\omega_{F_3}=O(x^{n-2k+2})$ since the first positive indicial root for $\Delta_{k-1}$ 
is $n-2k+2$. Using that $\delta_0L_k\omega_0=0$ and the form of $\delta_g$ we obtain 
\[\delta_g\omega_{F_3}=\delta_g\omega_{F_1}+x^{n-2k+2}\Big(\delta_0v^{(t)}-\frac{1}{n+2-2k}(B_k\omega_0+2\delta_0v^{(t)})\Big)+
O(x^{n-2k+3}).\]
By Proposition \ref{propr}, $\delta_g\omega_{F_1}=O_t(x^{n-2k+2})+O_n(x^{n-2k+4})$ and from the definition of $C_k$, 
a necessary condition to have $\delta_g\omega_F=O(x^{n-2k+3})$ is
\[(n-2k)\delta_0v^{(t)}=B_k\omega_0-(n-2k+2)C_k\omega_0.\]
Writing now $\delta_0v^{(t)}$ in terms of $B_k,C_k$ in \eqref{omegaF3} proves the uniqueness and the form of the expansion.
Now for the existence, one can take the form in Proposition \eqref{poisson}. Another way, which again is formal, is first to construct a polyhomogeneous $(k+1)$-form $\omega'_F$ such that
\[\left\{\begin{array}{l}
\Delta_{k+1}\omega'_F=O_t(x^{n-2k-1})+O_n(x^{n-2k+1})\\
\omega'_F=\frac{2(-1)^{k+1}}{n-2k}\log(x)d\omega_0 +\omega_0\wedge \frac{dx}{x}+O(x),
\end{array}\right.\] 
which can be done as in Proposition \ref{construQ} by using the indicial equations, 
and then to set $\omega_F:=\delta_g\omega'_F$. It is easy to see that this form is a polyhomogeneous solution of \eqref{condomegaF}.  
\qed\\

Since the exact solution in Proposition \ref{poisson} is coclosed, we deduce from Proposition \ref{formalequations} the 
\begin{cor}\label{hklocal}
The operator $G_k$ is a natural differential operator of order $n-2k+1$ which is given by 
\[G_k=(-1)^{k+1}\frac{B_k-2C_k}{n-2k}\]
and depends only on $h_0$ and the covariant derivatives of its curvature tensor. 
\end{cor}  

\subsubsection{The operator $Q_k$}\label{theoperatorqk} Following the ideas used above for $G_k$, we shall show how to construct $Q_k$
from a formal solution $\omega_{F_1}$. We start by  
\begin{defi}\label{B'k}
For $1\leq k<\ndemi-1$, define the operators $B'_{k-1}:C^\infty(M,\Lambda^{k-1}(M))\to C^{\infty}(M,\Lambda^{k-1}(M))$  
and $D'_{k-1}:C^\infty(M,\Lambda^{k-1}(M))\cap\ker d\to C^{\infty}(M,\Lambda^{k}(M))$ by
\begin{equation}
\begin{gathered}
B'_{k-1}\omega_0:=\Big(x^{-n+2k-2}i_{x\pl_x}\Delta_k\omega'_{F_1}\Big)|_{x=0},\\
D'_{k-1}\omega_0:=\Big(x^{-n+2k}i_{x\pl_x}d\omega'_{F_1}\Big)|_{x=0}
\end{gathered}
\end{equation}
where $\omega'_{F_1}$ is the form in \eqref{F2F1} such that $\Delta_{k}\omega'_{F_1}=O(x^{n-2k+2})$ and 
$\omega'_{F_1}=\omega_0\wedge \frac{dx}{x}+O(x^2)$. 
\end{defi}

Let us now set $\omega'_{F_2}:=\omega'_{F_1}+v^{(t)}x^{n-2k}$ for some arbitrary smooth form $v^{(t)}$ on $M$, 
we obtain
\[\begin{gathered}
\Delta_k\omega_{F_2}'=(-1)^{k+1}x^{n-2k+2}(B'_{k-1}\omega_0+2\delta_0v^{(t)})\wedge\frac{dx}{x}
+O_t(x^{n-2k+2})+O_n(x^{n-2k+3}).
\end{gathered}\]
so to solve away the $x^{n-2k+2}$ normal coefficient, we need to define
\begin{equation}\label{omegaF3bis}
\omega'_{F_3}:=\omega'_{F_2}+\frac{(-1)^{k+1}}{n-2k+2}x^{n-2k+2}\log(x)\Big(B_{k-1}'\omega_0+2\delta_0v^{(t)}\Big)\wedge
\frac{dx}{x}
\end{equation}
which satisfies $\Delta_{k}\omega'_{F_2}=O_t(x^{n-2k+2}\log(x))+O_n(x^{n-2k+3})$. Like for $G_k$, the term $v^{(t)}$
is arbitrary and so we have to impose an additional condition to fix this term (or at least to fix $\delta_0v^{(t)}$).
\begin{prop}\label{formalequations2}
Let $\omega_0\in C^\infty(M,\Lambda^{k-1}(M))$ be closed, then there is a polyhomogeneous $k$-form $\omega'_F$
which satisfies 
\[\left\{\begin{array}{l}
\Delta_k\omega'_{F}=O_t(x^{n-2k+1})+O_n(x^{n-2k+3})\\
d\omega'_{F}=O(x^{n-2k+1})\\
\omega'_F=\omega_0\wedge\frac{dx}{x}+O(x^2)\end{array}\right.,\]
which is unique modulo $O_t(x^{n-2k})+O_n(x^{n-2k+2})$ and has an expansion of the form
\begin{equation*}
\begin{split}
\omega'_F=&\,\sum_{j=1}^{\frac{n}{2}-k-1}x^{2j}\omega^{(t)}_{2j}+
\sum_{j=0}^{\frac{n}{2}-k+1}x^{2j}\omega^{(n)}_{2j}\wedge\frac{dx}{x}-x^{n-2k}\frac{1}{n-2k}D'_{k-1}\omega_0\\
&+\frac{(-1)^{k+1}}{n-2k+2}x^{n-2k+2}\log(x)\Big(B'_{k-1}\omega_0-\frac{2\delta_0D'_{k-1}\omega_0}{n-2k}\Big)\wedge\frac{dx}{x}.
\end{split}
\end{equation*}
\end{prop}
\textsl{Proof}: (i) Take $\omega'_F=\omega'_{F_3}$ defined in \eqref{omegaF3bis}, then $\Delta_k\omega'_F=
O_t(x^{n-2k+1})+O_n(x^{n-2k+3})$ by construction. Moreover, since $\omega_0$ is closed, one has
$d\omega'_F=O(x^2)$ and $\Delta_{k+1}d\omega'_F=O(x^{n-2k+1})$. Since the indicial roots for $\Delta_{k+1}$
in $[2,n-2k+1]$ are $n-2k-2$ in the $\Lambda^{k+1}_t$ part and $n-2k$ in the $\Lambda^{k+1}_n$ part, this implies that
$d\omega'_F=O_t(x^{n-2k-2})+O_n(x^{n-2k})$.
Then, using \eqref{domegaF1}, we obtain 
\begin{equation*}\begin{split}
d\omega'_F=&x^{n-2k}\Big(dv^{(t)}+\big((-1)^{k}(n-2k)v^{(t)}+d\omega^{(n)}_{n-2k}\big)\wedge\frac{dx}{x}\Big)+O(x^{n-2k+1})\\
=&\, x^{n-2k}\Big(dv^{(t)}+\big((-1)^{k}(n-2k)v^{(t)}+(-1)^kD'_{k-1}\omega_0\big)\wedge\frac{dx}{x}\Big)+O(x^{n-2k+1}).
\end{split}
\end{equation*}
So $d\omega'_F=O(x^{n-2k+1})$ if and only if $v^{(t)}=-D'_{k-1}\omega_0/(n-2k)$. 
\qed\\

The first corollary is 
\begin{cor}\label{cor1}
For $k<\ndemi$, the operator $Q_k$ is a natural differential operator of order $n-2k$ which is given by
\[Q_k=\frac{(-1)^k}{n-2k}\Big(B'_k\omega_0-\frac{\delta_0D'_k}{\ndemi-k-1}\Big)\] 
and it depends only on $h_0$ and the covariant derivatives of its curvature tensor.
\end{cor}

As a corollary of Propositions \ref{formalequations} and \ref{formalequations2}, we also have
\begin{cor}\label{cor2}
If $\omega_0$ is a  closed $k$-form on $M$, then there is a polyhomogeneous $k$-form $\omega_F$ on $\bar{X}$ such that
\[\left\{\begin{array}{l}
d\omega_{F}=O(x^{n-2k+1}))\\
\delta_g\omega_{F}=O(x^{n-2k+3})\\
\omega_F=\omega_0+O(x)\end{array}\right.\]
It is unique modulo $O_t(x^{n-2k+1})+O_n(x^{n-2k+2})$ and it has an expansion
\begin{equation*}\begin{split}
\omega_F=\,&\sum_{j=0}^{\frac{n}{2}-k-1}x^{2j}\omega^{(t)}_{2j}
+\sum_{j=1}^{\frac{n}{2}-k}x^{2j}\omega^{(n)}_{2j}\wedge\frac{dx}{x}-\frac{1}{n-2k}D_k\omega_0 x^{n-2k}\\
&+x^{n-2k+2}\log(x)\Big(\frac{(-1)^{k+1}}{n-2k}(B_k\omega_0-2C_k\omega_0)\wedge\frac{dx}{x}\Big).
\end{split}
\end{equation*}
\end{cor}
\textsl{Proof}: for the existence, take $\omega'_F$ in  Proposition \ref{formalequations2} ($\omega'_F$ is $k+1$ form now since $\omega_0\in \Lambda^{k}(M)$) and consider
$\omega_F:=(-1)^{k+1}/(2k-n)\delta_g\omega'_F$. It is easy to see that $\omega_F=\omega_0+O(x^2)$ and
that $\Delta_k\omega_F=O_t(x^{n-2k+1})+O_n(x^{n-2k+3})$.
Since $d\delta_g\omega'_F=-\delta_gd\omega'_F+O_t(x^{n-2k-1})+O_n(x^{n-2k+1})$,
we deduce that $d\omega_F=O_t(x^{n-2k-1})+O_n(x^{n-2k+1})$.
But from the Proposition \ref{formalequations}, $\omega_F=\omega_{F_1}+v^{(t)}x^{n-2k}+O(x^{n-2k+1})$ (note that $L_k\omega_0=0$ by Proposition \ref{propr2})
for some $k$-form $v^{(t)}$ on $M$ and so we conclude that
\begin{equation*}
\begin{split}
d\omega_F=&\,\sum_{2j=2}^{n-2k-2} d\omega^{(t)}_{2j}x^{2j}+
\sum_{2j=2}^{n-2k-2}x^{2j}((-1)^k2j\omega^{(t)}_{2j}+d\omega^{(n)}_{2j})\wedge\frac{dx}{x}\\
& +x^{n-2k}(d\omega_{n-2k}^{(n)}+(-1)^k(n-2k)v^{(t)})\wedge \frac{dx}{x}+x^{n-2k}dv^{(t)}+O(x^{n-2k+1})\\
=& \, O(x^{n-2k+1})
\end{split}
\end{equation*}
so $v^{(t)}$ has to be $(-1)^{k+1}d\omega_{n-2k}^{(n)}/(n-2k)$ to get $d\omega_F=O_t(x^{n-2k-1})+O_n(x^{n-2k+1})$.
But clearly this argument also implies that $d\omega_{F_1}=x^{n-2k}d\omega_{n-2k}^{(n)}\wedge\frac{dx}{x}$ and the expansion of $\omega_{F}$ is then a consequence of this fact together with the expansion \eqref{omegaF} in Proposition \ref{formalequations} and the definition of $D_k$.  
\qed\\

\textbf{Remark}: in Proposition \ref{formalequations}, \ref{formalequations2} and Corollary \ref{cor2}, we do not really need 
to take $\omega_0\in C^{\infty}(M,\Lambda(M))$. Indeed, for an $\omega_0$ in $L^2(\Lambda(M))$, the arguments
would work in a similar fashion except that the expansion in power of $x$ and $\log(x)$ have coefficients in some
$H^{-N}(\Lambda(M))$ with $N$ large enough, like we discussed in the proof of Proposition \ref{poisson}. 

\subsection{Factorizations}

\begin{prop}\label{propr2}
For any $k<\ndemi-1$, the following identities hold
\begin{equation}\label{relation}
\begin{split}
G_k=&(-1)^{k}\frac{\delta_{h_0}Q_k}{n-2k}\quad\,\,\,\textrm{on closed forms},\\
L_k=&\frac{(-1)^k}{(n-2k)}G_{k+1}d=-\frac{\delta_{h_0}Q_{k+1}d}{(n-2k)(n-2k-2)}.
\end{split}
\end{equation}
while for $k=\ndemi-1$ 
\begin{equation}\label{Ln}
L_{\ndemi-1}=\frac{1}{2}\delta_{h_0}d.
\end{equation}
\end{prop}
\textsl{Proof}: Let $\omega$ be the solution of \eqref{dx/x} with initial data 
$\omega_0$ closed. Then its first log term is $x^{n-2k+2}\log(x)Q_{k-1}\omega_0\wedge\frac{dx}{x}$ and thus the first normal log term
of $\delta_g\omega$ is $x^{n-2k+4}\log(x)(\delta_0Q_{k-1}\omega_0)\wedge\frac{dx}{x}$. But $\delta_g\omega$ is in 
$C^{n-2k+1}(\bar{X},\Lambda^{k-1}(\bar{X}))$ and is harmonic with leading behaviour at the boundary  
\[\delta_g\omega=(-1)^{k}(2k-n-2)\omega_0+O(x).\]
Thus,
the form $\delta_g\omega$ has for first normal log term $(-1)^k(2k-n-2)x^{n-2k+4}(G_{k-1}\omega_0)\wedge\frac{ dx}{x}$.

Since $\Delta_{k+1}d=d\Delta_{k}$ then $\omega':=d\omega$ is a solution (unique modulo ${\rm Ker}_{L^2}\Delta_{k+1}$) of $\Delta_{k+1}\omega'=0$ with
$\omega'|_{x=0}=d\omega_0$. 
But since the first log term in $d\omega$ is
\[(-1)^k(n-2k)x^{n-2k}\log(x)L_k\omega_0\wedge\frac{dx}{x}\]
 and since $L^2$ harmonic forms have no log term at this order, we get (\ref{relation}).

To compute $L_{\ndemi-1}$, we compute iteratively $\omega_{F_1}=\omega_0-x\frac{(-1)^{\ndemi}}{2}\delta_0\omega_0\wedge dx$, therefore  $\Delta_{\ndemi-1}\omega_{F_1}=x^2\delta_0d\omega_0+o(x^2)$. Which gives the result by \eqref{LkF1}.
\qed\\

\textbf{Remark}: Note that it implies that $L_k$ is zero on closed forms and $G_k$ has its range in co-closed forms.

\subsection{Conformal properties}\label{propoflk}

A priori our construction of $L_k,G_k,Q_k$ depends on the choice of geodesic boundary defining function $x$, i.e. 
on the choice of conformal representative in $[h_0]$.
In order to study the conformal properties of these operators, we need to compare the splittings of the differential forms associated to different conformal representatives.\\

A system of coordinates $y=(y_i)_{i=1,\dots,n}$ on $M$ near a point $p\in M$
give rise to a system of coordinates $(x,y)$ in $\bar{X}$ near the boundary point $p$
through the diffeomorphism $\psi:(x,y)\to \psi_x(y)$ where $\psi_t$ is the flow of the gradient
$\nabla^{x^2g}x$ of $x$ with respect to $x^2g$. Such system $(x,y)$ is called a system of 
\emph{geodesic normal coordinates associated to $h_0$}. 

\begin{lem}\label{changconf} 
Let $(x,y)$ and $(\hat{x},\hat{y})$ be two systems of geodesic normal coordinates associated respectively to $h_0$ and $\hat{h}_0=e^{2\varphi_0}h_0$. If $\hat{\omega}$ (resp. $\hat{\omega}\wedge d\hat{x}$) is a $k$-form tangential (resp. normal) in the coordinates $(\hat{x},\hat{y})$ with $\hat{\omega}|_{\hat{x}=0}=\omega_0$, then we have
$$\displaylines{\hat{\omega}=\omega_0+(-1)^{k+1}x^2(i_{\nabla\varphi_0}\omega_0)\wedge \frac{dx}{x}+O_t(x^2)+O_n(x^3),\cr
\hat{\omega}\wedge \frac{d\hat{x}}{\hat{x}}=\omega_0\wedge \frac{dx}{x}+\omega_0\wedge d\varphi_0+O_t(x)+O_n(x^2).}$$ 
\end{lem}
\textsl{Proof}: 
By the proof of Lemma 2.1 in \cite{Gu}, if $\hat{h}_0=e^{2\varphi_0}h_0$ is another conformal representative, a geodesic boundary defining function $\hat{x}$ associated to $\hat{h_0}$ satisifies $\hat{x}=e^{\varphi}x$ with $\varphi=\varphi_0+O(x^2)$ at least $C^{n-1}$ and $\hat{y}_i(x,y)=y_i+\frac{x^2}{2}\sum_{j=1}^nh^{ij}\pl_{y_j}\varphi_0+O(x^3)$. Hence $d\hat{y}_i=dy_i+x\sum_jh^{ij}\pl_{y_j}\varphi_0 dx$ and $d\hat{x}=xe^{\varphi_0} d\varphi_0+e^{\varphi_0} dx+O(x^2)$, which gives the relations above.\qed\\

This implies the following corollary:
\begin{cor}\label{lkgk}
Under a conformal change $\hat{h}_0=e^{2\varphi_0}h_0$, the associated operators $\hat{L}_k$, $\hat{H}_k$ and $\hat{Q}_k$ are given by 
\begin{equation}\label{changelaw}
\begin{gathered}
\hat{L}_k=e^{(2k-n)\varphi_0}L_k, \quad \hat{G}_k=e^{(2k-2-n)\varphi_0}\bigl(G_k+(-1)^ki_{\nabla\varphi_0}L_k\bigr)\\
\hat{Q}_{k}\omega_0=e^{\varphi_0(2k-n)}\Big(Q_k\omega_0+(n-2k)L_{k}(\varphi_0\omega_0)\Big)
\end{gathered}
\end{equation}
where $\omega_0\in C^\infty(M,\Lambda^k(M))$ is any closed form. Thus $L_k$ is conformally covariant and $G_k$ is conformally covariant on the kernel of $L_k$ (hence on closed forms).
\end{cor}
\textsl{Proof}: The solution $\omega$ in Proposition \ref{poisson} is unique up to ${\rm ker}_{L^2}(\Delta_k)$ which is composed of functions which are $O(x^{n-2k+2})$, so by Lemma \ref{changconf}, when we change $h_0$ to $\hat{h}_0$ the first $\log x$ term (i.e. the $x^{n-2k}\log x$ term) in the expansion of $\omega$ changes by a multiplication by $e^{(2k-n)\varphi_0}$. As for the $x^{n-2k+2}\log x$ term in the normal part, we have a similar
effect but the tangential $\hat{x}^{n-2k}\log \hat{x}$ term gives rise to a $x^{n-2k+2}\log x$ term which gives the term $i_{\nabla\varphi_0} L_k$.

Using Lemma \ref{changconf} in the expansion \eqref{omdx/x}, we obtain that the form $\omega$ solution of \eqref{dx/x} can be written  
\begin{equation*}
\begin{split}
\omega=&\,\,\omega_0\wedge\frac{d\hat{x}}{\hat{x}}-\omega_0\wedge d\varphi_0+
\sum_{j=1}^{n-2k}\hat{x}^j({\omega_j'}^{(t)}+{\omega_{j}'}^{(n)}\wedge \frac{d\hat{x}}{\hat{x}})\\
&+e^{-\varphi_0(n-2k+2)}\hat{x}^{n-2k+2}\log(\hat{x})(Q_{k-1}\omega_0)\wedge\frac{d\hat{x}}{\hat{x}}+O(\hat{x}^{n-2k+2}). 
\end{split}
\end{equation*}
Now since this is a solution of $\Delta_k\omega=0$ with leading behaviour $\omega_0\wedge(d\hat{x}/\hat{x}-d\varphi_0)$,
we can consider the Dirichlet problems \eqref{dirichletpb} and \eqref{dx/x} with 
the choice of boundary defining function $\hat{x}$, and by the uniqueness of their solution modulo $\ker_{L^2}(\Delta_k)$, 
we deduce that $\omega=\hat{\omega}_1+\hat{\omega}_2$ where $\hat{\omega_1}$ is the solution of \eqref{dirichletpb}
with initial data $-\omega_0\wedge d\varphi_0$, and $\hat{\omega}_2$
is the solution of \eqref{dx/x} with initial data $\omega_0$ and boundary defining function $\hat{x}$.
Consequently, one has $\hat{\omega}_2=\omega-\hat{\omega_1}$ and the $\hat{x}^{n-2k+2}\log \hat{x}$ normal term in $\hat{\omega}_2$ 
is given by 
\[\hat{Q}_{k-1}\omega_0\wedge \frac{d\hat{x}}{\hat{x}}=e^{-\varphi_0(n-2k+2)}Q_{k-1}\omega_0 \wedge \frac{d\hat{x}}{\hat{x}}+\hat{G}_k(\omega_0\wedge d\varphi_0)\wedge \frac{d\hat{x}}{\hat{x}}.\]
Now we use Corollary \ref{lkgk} and \eqref{relation} with $d\omega_0=0$ to see that
\begin{equation*}
\begin{split}
e^{\varphi_0(n-2k+2)}\hat{G}_k(\omega_0\wedge d\varphi_0)=
&(-1)^{k-1}G_kd(\varphi_0\omega_0)+(-1)^ki_{\nabla\varphi_0}L_kd(\varphi_0\omega_0)\\
=&(n-2k+2)L_{k-1}(\varphi_0\omega_0)
\end{split}
\end{equation*}
This ends the proof of the transformation law of $Q_{k-1}$ by conformal change.
\qed\\

\noindent\textbf{Remark}: while $Q_k$ on $\ker d$ is not conformally invariant (by Proposition \ref{lkgk}), the pairing  
$\cjg Q_ku,u\cjd_{L^2({\rm dvol}_{h_0})}$ for the metric $h_0$ is conformally invariant for $u\in \ker d$. Indeed, using \eqref{changelaw}, a 
conformal change of metric $\hat{h}_0=e^{2\varphi_0}h_0$ gives
\[\int_M \cjg \hat{Q}_ku,u\cjd_{\hat{h}_0}{\rm dvol}_{h_0}=\int_M\cjg Q_ku,u\cjd_{h_0}+\frac{\cjg\delta_{0}Q_{k+1}d(\varphi_0u),u\cjd_{h_0}}{2k+2-n}{\rm dvol}_{h_0}\]
which by integration by part and $du=0$ gives the $\cjg \hat{Q}_ku,u\cjd_{L^2({\rm dvol}_{\hat{h}_0})}=
\cjg Q_ku,u\cjd_{L^2({\rm dvol}_{h_0})}$. Of course, when we restrict this form to exact forms, this is given by
\[\cjg Q_kdu,du\cjd=\cjg L_{k-1}u,u\cjd\]
which is real and conformally invariant.

\subsection{Analytical properties}

\subsubsection{Principal parts}

\begin{prop}\label{symb}
For any $k<\ndemi$ we have
$$\displaylines{Q_k=\frac{(-1)^{\frac{n}{2}+k+1}(n-2k)(\Delta_0)^{\frac{n}{2}-k}}{2^{n-2k}[(\frac{n}{2}-k)!]^2}+\mbox{\sl lower order terms in }\pl^j_{y_i}\cr
L_k=\frac{(-1)^{\frac{n}{2}+k+1}(n-2k)(\delta_0d)^{\frac{n}{2}-k}}{2^{n-2k}[(\frac{n}{2}-k)!]^2}+\mbox{\sl lower order terms in }\pl^j_{y_i}\cr
G_k=\frac{(-1)^{\frac{n}{2}+1}(\delta_0d)^{\frac{n}{2}-k}\delta_0}{2^{n-2k}[(\frac{n}{2}-k)!]^2}+\mbox{\sl lower order terms in }\pl^j_{y_i}}$$
\end{prop}
\noindent\textsl{Proof}: We first precise the computation of $\omega_{F_1}$ which solves \eqref{P'}. 
By Lemma \ref{formlapl}, $\omega_{F_1}$ has the form $\omega_{F_1}=\sum_{i=0}^{\frac{n}{2}-k-1}x^{2i}\omega_{2i}^{(t)}+\sum_{i=1}^{\frac{n}{2}-k}x^{2i}\omega_{2i}^{(n)}\wedge\frac{dx}{x}$, where the $\omega_i^{(*)}$ are images of $\omega_0$ by differential opertors on $M$. We compute the principal part of these operators by recurrence.

The decomposition \eqref{lapgrad} of $P$ and the identity $\Delta_k\omega_{F_1}=O_t(x^{n-2k})+O_n(t^{n-2k+1})$ give
\[\displaylines{
\hfill\sum_{i=1}^{\frac{n}{2}-k}x^{2i}\Bigr(-4i(k+i-\frac{n}{2}-1)\omega^{(n)}_{2i}+\sum_{j=1}^{i-1}Q'_j\omega^{(t)}_{2i-2j-2}+\bigl(R'_{j}+(k+2i-2j-1)P'_{j}\bigr)\omega^{(n)}_{2i-2j}\Bigr)\wedge\frac{dx}{x}\cr
\hfill+\sum_{i=0}^{\frac{n}{2}-k-1}x^{2i}\Bigl(-4i(k+i-\frac{n}{2})\omega^{(t)}_{2i}+\sum_{j=1}^{i}\bigl(R_{j}+(k+2i-2j-2)P_{j}\bigr)\omega^{(t)}_{2i-2j}+\sum_{j=1}^{i-1}Q_j\omega^{(n)}_{2i-2j}\Bigr)=0}\]
This determines uniquely the $\omega^{(*)}_i$.

Let us write ${\rm LOT}$ for lower order term operators on $M$.
Then we get 
\[\omega^{(n)}_2=\frac{(-1)^{k+1}}{2k-n}\delta_0\omega_0,\quad \omega^{(t)}_2=\Bigl(\frac{d\delta_0}{2(2k-n)}+\frac{\delta_0 d}{2(2k+2-n)}+{\rm LOT}\Bigr)\omega_0\] 
and given the order in $\pl_{y_i}$ of the $R_i,R'_i,\bar{R}_i,\bar{R}'_i,Q_i$ and $Q'_i$, we have
\[\begin{gathered}
\omega^{(t)}_{2i}=\frac{1}{2i(2k+2i-n)}\Bigl(2(-1)^{k+1}d\omega^{(n)}_{2i}+\Delta_0\omega^{(t)}_{2i-2}\Bigr)+
{\rm LOT}(\omega_0)\cr
\omega^{(n)}_{2i+2}=\frac{1}{2(i+1)(2k+2i-n)}\Bigl(2(-1)^{k+1}\delta_0\omega^{(t)}_{2i}+\Delta_0\omega^{(n)}_{2i}\Bigr)+
{\rm LOT}(\omega_0)
\end{gathered}\]
So we have
\[\begin{gathered}
\omega^{(t)}_{2i}=\Bigl(a_{2i}(\delta_0d)^i+b_{2i}(d\delta_0)^i+
{\rm LOT}\Bigr)\omega_0\cr
\omega^{(n)}_{2i+2}=\Bigl(a_{2i+1}(\delta_0d)^l\delta_0+{\rm LOT}\Bigr)\omega_0\end{gathered}\]
where the sequences $(a_i)$ and $(b_{2i})$ satisfy the relations
\[\begin{gathered}
a_{2i}=\frac{a_{2i-2}}{2i(2k+2i-n)}, \quad\quad\quad a_{2i+1}=\frac{2(-1)^{k+1}b_{2i}}{2(i+1)(2k+2i-n)}+\frac{a_{2i-1}}{2(i+1)(2k+2i-n)}\\
b_{2i}=\frac{2(-1)^{k+1}a_{2i-1}}{2i(2k+2i-n)}+\frac{b_{2i-2}}{2i(2k+2i-n)}
\end{gathered}\]
and $a_1=\frac{(-1)^{k+1}}{2k-n}$, $a_2=\frac{1}{2(2k+2-n)}$, $b_2=\frac{1}{2(2k-n)}$. By uniqueness of the solution of this equation we find
\[\begin{gathered} a_{2i}=\frac{1}{2^ii!\prod_{j=1}^i(2k+2j-n)},\quad\quad a_{2i+1}=\frac{(-1)^{k+1}}{2^ii!\prod_{j=0}^i(2k+2j-n)},\\ b_{2i}=\frac{1}{2^ii!\prod_{j=0}^{i-1}(2k+2j-n)}\end{gathered}\]
for all $i\leq\frac{n}{2}-k-1$.
We infer the equality
\begin{equation}
\begin{split}
\Delta_k\omega_{F_1}=&x^{n-2k}
\Bigl(a_{n-2k-2}(\delta_0d)^{\frac{n}{2}-k}+(b_{n-2k-2}+2(-1)^{k+1}a_{n-2k-1})(d\delta_0)^{\frac{n}{2}-k}+
{\rm LOT}\Bigr)\omega_0\\
 &+x^{2k-n+2}\Bigl(a_{n-2k-1}(\delta_0d)^{\frac{n}{2}-k}\delta_0+{\rm LOT}\Bigr)\omega_0\wedge\frac{dx}{x}+o(x^{n-2k+1})\\
=&x^{n-2k}\Bigl(\frac{(\delta_0d)^{\frac{n}{2}-k}}{2^{\frac{n}{2}-k-1}(\frac{n}{2}-k-1)!\prod_{j=1}^{\frac{n}{2}-k-1}(2k+2j-n)}+{\rm LOT}\Bigr)\omega_0\\
 &+x^{2k-n+2}\Bigl(\frac{(-1)^{k+1}(\delta_0d)^{\frac{n}{2}-k}\delta_0}{2^{\frac{n}{2}-k-1}(\frac{n}{2}-k-1)!
\prod_{j=0}^{\frac{n}{2}-k-1}(2k+2j-n)}+{\rm LOT}\Bigr)\omega_0\wedge\frac{ dx}{x}\\
 &+o(x^{n-2k+1})
\end{split}
\end{equation}
so we have 
\[\begin{gathered}
L_k=\frac{-(\delta_0d)^{\frac{n}{2}-k}}{2^{\frac{n}{2}-k-1}(\frac{n}{2}-k-1)!\prod_{j=0}^{\frac{n}{2}-k-1}(2k+2j-n)}+
{\rm LOT}\\
B_k=\frac{(\delta_0d)^{\frac{n}{2}-k}\delta_0}{2^{\frac{n}{2}-k-1}(\frac{n}{2}-k-1)!\prod_{j=0}^{\frac{n}{2}-k-1}(2k+2j-n)}+{\rm LOT}\end{gathered}\]
Note also that $\delta_g$ is of order $1$ so $C_k$ has no contribution to the principal part of $G_k$ and we get
\[G_k=\frac{(-1)^{k+1}(\delta_0d)^{\frac{n}{2}-k}\delta_0}{2^{\frac{n}{2}-k}(\frac{n}{2}-k)!\prod_{j=0}^{\frac{n}{2}-k-1}(2k+2j-n)}+{\rm LOT}.\]
The proof is the same (and even easier) for $Q_k$. We could have deduced the principal parts of $L_k$ and $G_k$ from the one of $Q_k$, but a slight generalization of the proof above will allow to compute the principal part of the non-critical $L_k^l$ in the next section.
\qed\\

We finally that the operators $L_k$ and $Q_k$ are symmetric on $C^\infty(M,\Lambda(M))$:
\begin{prop}\label{symmetry}
For $k\leq \ndemi-1$, the operators $L_k$ are symmetric on $C^{\infty}(M,\Lambda^k(M))$ while for $k<\ndemi-1$, the operators $Q_k$ are symmetric on $C^{\infty}(M,\Lambda^k(M))\cap \ker d$.
\end{prop}
\textsl{Proof}: The proof for $L_k$ is done in Proposition \ref{selfad} which covers the non-critical cases. 
The proof for $Q_k$ is quite similar, we let $\omega_0,\omega_0'$ be two closed $k$-forms on $M$ and 
$\omega,\omega'$ the forms constructed in the proof of Proposition \ref{construQ} with respective initial conditions
$\omega_0$ and $\omega_0'$. 
Then integration by part and the fact that $d\omega=d\omega'=0$ gives
\[\begin{gathered}
0=\int_{x\geq \eps}(\cjg \Delta_k\omega,\omega'\cjd_g-\cjg \Delta_k\omega',\omega\cjd_g)
{\rm dvol}_g=\\
\int_{x=\eps}\Big(\cjg i_{x\pl_x}\omega,\delta_g\omega'\cjd_{h_x}-\cjg i_{x\pl_x}\omega',\delta_g\omega\cjd_{h_x}\Big)
x^{-n}{\rm dvol}_{h_x}.
\end{gathered}\]
But  a straightforward analysis and the fact that $L_k(\omega_0)=L_k(\omega'_k)=0$ give that the second line has an expansion of the form 
\[\begin{gathered}
a_{-2\ell}\eps^{-2\ell}+\dots+ a_{-2}\eps^{-2}
+L\log(\eps)+O(1)\\
\textrm{ with }L:=(-1)^{k+1} (2k-n)\Big(\cjg Q_k\omega_0,\omega_0'\cjd_{L^2({\rm dvol}_{h_0})}-\cjg \omega_0,Q_k\omega_0'\cjd_{L^2({\rm dvol}_{h_0})}\Big)
\end{gathered}\] 
This achieves the proof.\qed  

\subsection{Branson $Q$-curvature}\label{qcurvature}

We finally conclude this section by the observation that $Q_0$ is the $Q$-curvature of Branson.
\begin{prop}\label{qbranson}
The operator $Q_0$ of Definition \ref{definQ} satisfies
\[Q_01=\frac{n(-1)^{\ndemi+1}}{2^{n-1}\ndemi!(\ndemi-1)!}Q\] 
where $Q$ is Branson $Q$-curvature defined in \cite{B}. 
\end{prop}
\textsl{Proof}: Let $(X,g)$ Poincar\'e-Einstein with conformal infinity $(M,[h_0])$. In \cite{FG3}, Fefferman and Graham showed that the $Q$-curvature of Branson is the function $Q$ on $M$ such that if $U\in C^\infty(X)$ is the function solution of
\[\left\{\begin{array}{l}
\Delta_gU=n\\
U=\log(x)+A+x^nB\log(x) \textrm{ with }A,B\in C^{\infty}(\bar{X})\\
A|_{x=0}=0
\end{array}\right.\]
then $B|_{x=0}=(-1)^{\ndemi+1}(2^{n-1}\ndemi!(\ndemi-1)!)^{-1}Q$. Consider $dU$, clearly it is a harmonic $1$-form
and it is given by
\[dU=\frac{dx}{x}+dA+nx^{n}B\log(x)\frac{dx}{x} +O(x^{n})\]
and by uniqueness of the solution in Proposition \ref{construQ} and the decay  
of $L^2$ harmonic $1$-forms (of order $x^n$), we deduce that $Q_01=nB|_{x=0}$, this proves the claim (note that the log term in the development of $\Delta_k$ does not interfer since it acts trivially on normal zero forms).
\qed


\section{The non-critical case}
Let $(X,g)$ be a Poincar\'e-Einstein manfiold with conformal infinity $(M,[h_0])$. 
We assume $k\leq (n+1)/2$ and $n$ may be odd or even in this section, and we let $\ell$ be an integer
in $[1,\ndemi-k]$ in general, and $\ell\in \nn$ if $n$ is odd and $(X,g)$ is an even 
Poincar\'e-Einstein manifold. 
We want to construct the operators $L_k^\ell$ of \cite{BG} by solving the following equation
\begin{equation}\label{noncrit}
\begin{gathered}
\Big(\Delta_k-(\ndemi-k+\ell)(\ndemi-k-\ell)\Big)\omega=O_t(x^{\ndemi-k+\ell})+O_n(x^{\ndemi-k+\ell+1})\\
\textrm{ with } \omega=x^{\ndemi-k-\ell}\omega_0+o(x^{\ndemi-k-\ell}) \textrm{ as }x\to 0.
\end{gathered}
\end{equation}
where $O_n,O_t$ are defined in the proof of Proposition \ref{poisson} and where $\omega_0\in C^{\infty}(M,\Lambda^k(M))$. 
This can be done essentially like in the critical case, using the indicial equations of Subsection \ref{indeq}.
Indeed, the indicial roots of $\Delta_k-(\ndemi-k+\ell)(\ndemi-k-\ell)$ can be computed rather easily, these are 
\[\begin{gathered}
\ndemi -k-\ell \quad\textrm{and} \quad\ndemi -k +\ell \quad \textrm{in the }\Lambda^k_t \textrm{ component }\\ 
\ndemi-k-\sqrt{\ell^2+n+1-2k}\quad  \textrm{and}\quad \ndemi-k+\sqrt{\ell^2+n+1-2k}\quad \textrm{in the }\Lambda_n^k \textrm{ component}.  
\end{gathered}\] 
Since there is no indicial roots in $(\ndemi-k-\ell,\ndemi-k+\ell)$, we obtain 
\begin{lem}\label{casnoncrit}
For $\omega_0\in C^{\infty}(M,\Lambda^{k}(M))$ fixed, there exists a series
\begin{equation}\label{omegaf1}
\omega_{F_1}=x^{\ndemi-k-\ell}\Big(\sum_{2j=0}^{2l-2}x^{2j}\omega^{(t)}_{2j}+
\sum_{2j=2}^{2l}x^{2j}(\omega_{2j}^{(n)}\wedge \frac{dx}{x})\Big)
\end{equation}
such that $\omega_0^{(t)}=\omega_0$ and
\begin{equation}\label{deltaomega}
\Big(\Delta_k-(\ndemi-k+\ell)(\ndemi-k-\ell)\Big)\omega_{F_1}=O_t(x^{\ndemi-k+\ell})+O_n(x^{\ndemi-k+\ell+2})
\end{equation} 
where the forms $\omega^{(.)}_j$ on $M$ are uniquely determined by $\omega_0$ and the expansion of $\Delta_k$
in powers of $x$ given by Lemma \ref{formlapl}. 
\end{lem}
Note that the condition $\ell\leq \ndemi$ insures that
that the first $\log(x)$ coefficient coming from the metric does not show up in \eqref{noncrit}.
Since $(\ndemi-k+\ell)$ is an indicial root in the $\Lambda^k_t$ component, we can then define 
\begin{equation}\label{omegan-k}
\begin{gathered}
\omega_{F_2}=\omega_{F_1}+x^{\ndemi-k+\ell}\log(x)\omega_{n-k+\ell,1}^{(t)}, \\
\textrm{ with }\omega_{n-k-\ell,1}^{(t)}=\frac{1}{2\ell}\Big[x^{-\ndemi+k-\ell}\Big(\Delta_k-(\ndemi-k+\ell)(\ndemi-k-\ell)\Big)\omega_{F_1}\Big]_{|_{x=0}}
\end{gathered}\end{equation}
which satisfies 
\begin{equation}\label{propomegaf2}
\Big(\Delta_k-(\ndemi-k+\ell)(\ndemi-k-\ell)\Big)\omega_{F_2}=O_t(x^{\ndemi-k+\ell+1})+
O_n(x^{\ndemi-k+\ell+2}\log x).
\end{equation}

\noindent\textbf{Remark}: we could continue the construction to get a solution $\omega$ of 
\[(\Delta_k-(\ndemi-k-\ell)(\ndemi-k+\ell))\omega=O(x^{\infty})\] 
and even an exact solution (with no $O(x^{\infty})$) using the resolvent of $\Delta_k$. However, 
since the mapping properties of $(\Delta_k-(\ndemi-k-\ell)(\ndemi-k+l))^{-1}$ is not really
available in the literature when $\ell\not=\ndemi-k$, we do not discuss this case further.\\

Like we did for $L_k$, we can then define an operator on $M$ as follows: 
\begin{defi}\label{Lkl} 
For $k\leq (n+1)/2$, we let $\ell$ be an integer in $[1,\ndemi]$ if $n$ is even and in $\nn$ if $n$ is odd.
The operator $L_{k}^\ell: C^{\infty}(M,\Lambda^{k}(M))\to C^{\infty}(M,\Lambda^{k}(M))$ is defined by
$L_k^\ell\omega_0:=\omega_{n-k-\ell,1}^{(t)}$ where $\omega_{n-k-\ell,1}^{(t)}$ is 
given in \eqref{omegan-k}. 
\end{defi}

\noindent\textbf{Remark}: clearly, we have $L_k^{\ndemi-k}=L_k$ when $n$ is even.
 
\begin{lem}\label{deltaom}
The form $\omega_{F_1}$ of \eqref{omegaf1} satisfies $\delta_{g}\omega_{F_1}=O(x^{\ndemi-k+\ell+2})$.
\end{lem}
\textsl{Proof}: by \eqref{deltaomega} and $\delta_g\Delta_k=\Delta_{k-1}\delta_g$, 
the form $\delta_g\omega_{F_1}$ solves 
\begin{equation}\label{approx}
\Big(\Delta_{k-1}- (\ndemi-k+\ell)(\ndemi-k-\ell)\Big)\delta_g\omega_{F_1}=
O(x^{\ndemi-k+\ell+2})
\end{equation}
and with $\delta_g\omega_{F_1}=O(x^{\ndemi-k-\ell+2})$. The Taylor series $T$ of $x^{-\ndemi+k+\ell}\delta_g\omega_{F_1}$ 
to order $O(x^{2\ell+2})$ is such that $x^{\ndemi-k-\ell}T$ 
solves \eqref{approx}, moreover $T$ is even by Lemma \ref{formlapl}.
A short computation shows that there is no
indicial roots of $(\Delta_{k-1}-(\ndemi-k+\ell)(\ndemi-k+\ell))$ in the interval 
$[\ndemi-k-\ell+2,\ndemi-k+\ell+1]$ except when $2k=n+1$ where $\ndemi-k+\ell+1$ is a root in the $\Lambda^t$ component, 
this implies that the Taylor series of $\delta_g\omega_{F_1}$ vanishes to order 
$O(x^{\ndemi-k+\ell+1})$ except maybe when $2k=n+1$. However in the last case, 
by parity of $T$, we see that there is no $\ndemi-k+\ell+1$ term in the expansion 
of $\delta_g\omega_{F_1}$, this ends the proof.
\qed\\

By an obvious integration by part, we have the
\begin{prop}\label{selfad} 
The operators $L_k^\ell$ are symmetric on $C^{\infty}(M,\Lambda^k(M))$.
\end{prop}
\textsl{Proof}:  Consider $\omega^1_{F_2}$ and $\omega_{F_2}^2$ like in \eqref{omegan-k} with respective boundary
values $\omega_0^1$ and $\omega_0^2$, they are well defined form in some collar neighbourhood $X_1:=(0,\eps_0)_x\x M$ 
of $M$ in $\bar{X}$. 
Let $\varphi\in C_0^{\infty}((-\eps_0,\eps_0))$ be a cut-off function which equals $1$ near $0$ and $\til{\omega}^i:=\varphi(x)\omega^i_{F_2}$
for $i=1,2$. Then using Lemma \ref{deltaom} we have $\delta_g\til{\omega}^i=O(x^{\ndemi-k+\ell+1})$, but since
$i_{x\pl_x}\til{\omega}^i=O(x^{\ndemi-k-\ell+2})$, the Green formula gives for small $\eps>0$
\[\begin{gathered}
\int_{x\geq \eps}(\cjg \Delta_k\til{\omega}^1,\til{\omega}^2\cjd_g-\cjg \Delta_k\til{\omega}^2,\til{\omega}^1\cjd_g)
{\rm dvol}_g=
(-1)^n\int_{x=\eps}(\star_gd\til{\omega}^1)\wedge\til{\omega}^2-(\star_g\til{\omega}^2)\wedge\til{\omega}^1
+O(\eps).
\end{gathered}\]
But the first line is a $O(1)$ as $\eps\to 0$ by \eqref{propomegaf2}, and a straightforward analysis
gives that the second line has an expansion of the form 
\[\begin{gathered}
a_{-2\ell-1}\eps^{-2\ell-1}+\dots+ a_{-1}\eps^{-1}
+L\log(\eps)+O(1)\\
\textrm{ with }L:=(-1)^n(\ndemi-k+\ell)\int_M (\star_{0}L_k^\ell\omega_0^1)\wedge\omega_0^2-(\star_{0}L_k^\ell\omega_0^2)\wedge \omega_0^1
\end{gathered}\]  
and this implies $L=0$ by comparing the $\log(\eps)$ terms.
\qed\\

\begin{lem}
  We have $L_k^l=\frac{(-1)^{l+1}l}{2^{2l-1}(l!)^2}\Bigl[(\delta_0d)^l+\frac{n-2k-2l}{n-2k+2l}(d\delta_0)^l\Bigr]+{\rm LOT}$.
\end{lem}
\textsl{Proof}:
We set $T$ such that $\omega_{F_1}=x^{\frac{n}{2}-k-l}T$, $\lambda=(\frac{n}{2}-k+l)(\frac{n}{2}-k-l)$ and $P=x^{k+l-\frac{n}{2}}(\Delta-\lambda)x^{\frac{n}{2}-k-l}$.

Then we have $T=\sum_{i=0}^{l-1}x^{2i}\omega_{2i}^{(t)}+\sum_{i=1}^{l}x^{2i}\omega_{2i}^{(n)}\wedge\frac{dx}{x}$ and $P$ admits the same decomposition than $\Delta_k$ in Lemma \ref{formlapl} but with indicial operator equal to
$$
\begin{pmatrix}
  2lx\pl_x-(x\pl_x)^2&2(-1)^{k+1}d\\
0&-(x\pl_x)^2+2(l+1)x\pl_x+n-2k-2l
\end{pmatrix}$$
The equation $PT=O_t(x^{2l})+O_n(x^{2l+1})$ gives then
$$\displaylines{\hfill\omega_{2i+2}^{(n)}=\bigl(a_{2i+1}(\delta_0d)^i\delta_0+{\rm LOT}\bigr)\omega_0\hfill\omega_{2i}^{(t)}=\bigl(a_{2i}(\delta_0d)^i+b_{2i}(d\delta_0)^i+{\rm LOT}\bigr)\omega_0\hfill}$$
with
$$
\begin{gathered}
  a_1=\frac{(-1)^k}{\frac{n}{2}-k+l},\quad\quad a_2=\frac{-1}{4(l-1)},\quad\quad b_2=\frac{-(\frac{n}{2}-k+l-2)}{4(l-1)(\frac{n}{2}-k+l)}
\end{gathered}$$
and
$$
\begin{gathered}
  a_{2i+2}=\frac{a_{2i}}{4(i+1)(i+1-l)},\quad\quad a_{2i+1}=\frac{2(-1)^{k+1}b_{2i}+a_{2i-1}}{4(i+1)(i-l)+2k-n+2l},\\ b_{2i+2}=\frac{b_{2i}+2(-1)^{k+1}a_{2i+1}}{4(i+1)(i+1-l)}.
\end{gathered}$$
The solutions of these equations are
$$
\begin{gathered}
  b_{2i}=\frac{(-1)^i(n-2k+2l-4i)(l-i-1)!}{4^ii!(l-1)!(n-2k+2l)}\quad\quad a_{2i+1}=\frac{(-1)^{k+i}(l-i-1)!}{2^{2i-1}i!(l-1)!(n-2k+2l)}\\ a_{2i}=\frac{(-1)^i(l-i-1)!}{4^ii!(l-1)!}
\end{gathered}$$
Since the equation \eqref{omegan-k} reads
$$L_k^l=\Bigl[\frac{x^{-2l}}{2l}i_{x\pl_x}\bigl(\frac{dx}{x}\wedge PT\bigr)\Bigr]_{|x=0}$$
we get the result.
\qed\\
\section{Relation with Branson-Gover operators}

First we recall a few fact on the ambient metric of Fefferman-Graham, 
see \cite{FG,FG2} for details.
If $(M,[h_0])$ is a compact manifold equipped with a conformal class, 
we call 
\[\mc{Q}=\{t^2h_0(m); t>0,m\in M\}\subset S^2T^*M\] 
the conformal bundle, it is identified with $(0,\infty)_t\x M$. Let $\til{\mc{Q}}=(-1,1)\x \mc{Q}$ be the ambient space with the inclusion $\iota:\mc{Q}\to \til{Q}$ defined by $z\to (0,z)$. There are
dilations $\delta_s: (t,m)\to (st,m)$ of $\mc{Q}$ which extends naturally to $\til{\mc{Q}}$.
The functions on $\mc{Q}$ which are $w$-homogeneous in the sense
\[f(st,m)=s^wf(t,m)\]
are the section of a bundle denoted $E[w]$, they extend naturally on $\til{\mc{Q}}$.
We denote by $\til{h}$ the ambient metric of Fefferman-Graham \cite{FG} on $\til{\mc{Q}}$.
This is a smooth Lorentzian metric on $\mc{Q}$ such that 
\[\begin{array}{ll}
(1) &\delta_s^*\til{h}=s^2\til{h}, \forall s>0,\\
(2)&\iota^*\til{h} \textrm{ is the tautological tensor }t^2h_0\textrm{ on }\mc{Q},\\
(3^*)& \Ric(\til{h})\textrm{ vanishes to infinite order at }\mc{Q} \textrm{ if }n \textrm{ is odd},\\
(3^{**})& \Ric(\til{h})\textrm{ vanishes to order }\ndemi-1\textrm{ at }\mc{Q}\textrm{ if } n\textrm{ is even}.
\end{array}\]
We let $T$ be the vector field which generates the dilations $\delta_s$, and let 
\[Q=\til{h}(T,T), \quad \rho:= -t^{-2}Q/2, \quad x=\sqrt{2\rho}, \quad u=xt\]
so that $Q$ is homegeneous of degree $2$ with respect to $\delta_s$, $u$  and $t$ are homogeneous of degree $1$
and $x$ of degree $0$, moreover $Q,\rho$ are smooth defining function of $\mc{Q}$, $x,u$ are defining function 
of $\mc{Q}$ in $\{Q\leq 0\}$ for some finer smooth structure on $\{Q\leq 0\}$.
Let us define $\mc{C}:=\{Q=-1, \rho<\eps\}$ for some small fixed $\eps$, then $\mc{C}$ can be identified
with a collar $(0,\eps)_\rho\x M$ and there is a system of coordinates $(u,m)\in (0,1]\x \mc{C}$
that covers the part $\{0>Q\leq -1, \eps>\rho>0\}$ which is a neighbourhood
of the cone $\mc{Q}$ near $t=\infty$. The metric $\til{h}$ has the model form (see \cite{FG})
in this neighbourhood 
\[\til{h}=-du^2+u^2g\]
where $g=(dx^2+h_x)/x^2$ is a Poincar\'e-Einstein metric on the collar $\mc{C}$.
  
The space $\mc{T}^{k}[s]$ is the space of $k$-form tractors which are homogeneous of degree $s$, i.e. these are 
restrictions to the null cone $\mc{Q}$ of $k$-forms on $\til{\mc{Q}}$ and such that $\til{\nabla}_TF=sF$
where $T=t\pl_t=u\pl_u$ is the generator of dilations in the cone fibers, $\til{\nabla}$ is the Levi-Civita connection 
on $\til{\mc{Q}}$. Since $\til{\nabla}_T *=*$ for $*=T,\pl_x, \pl_{m_i}$, we have $\mc{L}_T=\til{\nabla}_T+k$, on $\mc{T}^k[s]$, where $\mc{L}$ denotes Lie derivative.
The bundle  $\mc{E}^k[s]$ is the bundle that consists of 
the s-homogeneous $k$ forms on $M$, in the sense that they are the sections 
of $\Lambda^kT^*M\otimes E[s]$ and thus satisfy $\mc{L}_T\omega=s\omega$. We can view $\mc{E}^k[s]$
as a subspace of $\mc{T}^k[s-k]$.
We let $\mc{G}_k[s]$ be the subundle of $\mc{T}^{k}[s+k-n]$ consisting of forms which are annihilated by the interior product 
$i_T$. It has a conformally invariant projection onto $\mc{E}^k[s+2k-n]$ denoted by $q^k$, this is given for instance by 
$ i_{\pl_\rho}d\rho\wedge $.  

If $\til{\Delta}$ is the ambient Laplacian on $\til{\mc{Q}}$ associated to $\til{h}$,
if $\omega_0\in\mc{E}^k[k+\ell-\ndemi]$ and 
$\til{\omega}_0$ is an homogeneous extension of $\omega_0$ to $\til{\mc{Q}}$, then
it is proved in \cite[Prop. 4.3]{BG} that the operator defined by the formula
\begin{equation}\label{defL}
\mathbf{L}_k^\ell\omega_0=\Big[\iota_{T}\Big(\til{d}(n+2\til{\nabla}_T-2)+\demi \til{d}Q\wedge\til{\Delta}\Big)\til{\Delta}^\ell\til{\omega}_0\Big]_{|_{\mc{Q}}}=\Big[\iota_{T}\til{d}(n+2\til{\nabla}_T-2)\til{\Delta}^\ell\til{\omega}_0\Big]_{|_{\mc{Q}}}
\end{equation}
can be viewed as a conformally invariant operator $\mathbf{L}_k^\ell:\mc{E}^k[k+\ell-\ndemi]\to \mc{G}_k[\ndemi-k-\ell]$. 
Here $\til{d}$ denotes the exterior differential on $\til{\mc{Q}}$.
They also define the operators (see Proposition 4.4 and Theorem 4.5 in \cite{BG})  
\[L_k^{{\rm BG},\ell}:=q^k\mathbf{L}_k^\ell: \mc{E}^k[k+\ell-\ndemi]\to \mc{E}^k[k-\ndemi-\ell],\]  
\begin{equation}\label{Gkbrangov}
G^{\rm BG}_k:=q^{k-1}i_Y\mathbf{L}_k^{\ndemi-k}: \mc{E}^k[0]\to \mc{E}^{k-1}[2k-2-n]
\end{equation}
where $Y=-\frac{\pl_\rho}{t^2}$ is a vector field dual to $\til{d}t/t$ via $\til{h}$, it satisfies in particular $\til{d}Q(Y)=2$.
Finally the operator $Q^{\rm BG}_k$ acting on a closed $k$-form $\omega_0$ is defined as follows
\[Q_k^{\rm BG}\omega_0:=-2(\ndemi-k+1)q^k\Big[i_Yi_T\til{\Delta}^{\ndemi-k}(\frac{dQ}{2}\wedge \frac{\til{d}t}{t}\wedge\til{\omega}_0)\Big]|_{\mc{Q}}.\]
where $\til{\omega}_0$ is any homogeneous extension of $\omega_0$ to $\til{\mc{Q}}$.

We now prove a Lemma which is essentially the same than the proof for functions in \cite{GJMS}.
\begin{lem}\label{lem1}
Let $\omega\in\mc{T}^{k'}[-\alpha]$ and $j\in\mathbb{N}$, then we have
$$\til{\Delta}(Q^j\omega)=4j(\alpha-\frac{n}{2}-j)Q^{j-1}\omega+Q^j\til{\Delta}\omega.$$
\end{lem}
\textsl{Proof}:
Using $\til{\nabla} Q=2T$, we have $[\til{\Delta},Q]=-2(2\til{\nabla}_T+n+2)$ and so we can compute
\begin{equation*}
\begin{split}
\til{\Delta} (Q^{j}L\omega_0)=&\,\sum_{m=0}^{j-1}Q^m[\til{\Delta},Q]Q^{j-1-m}\omega+Q^j\til{\Delta}\omega\\
=&\,-2Q^{j-1}\sum_{m=0}^{j-1}\bigl(2(-2m+2j-2-\alpha)+n+2\bigr)\omega+Q^j\til{\Delta}\omega\\
=&\,4Q^{j-1}j(\alpha-\frac{n}{2}-j)\omega +Q^j\til{\Delta}\omega
\end{split}
\end{equation*}
which achieves the proof. \qed\\ 

As a consequence, and using Lemma \ref{lem1}, we get the 
\begin{theo}\label{coincide}
\textbf{(i)} Let $L_k^\ell, L_k$ and $G_k$ be the operators of Definition \ref{Lkl} and \ref{LkHk}, and let 
$c^{\ell}_k:=(-4)^\ell(\ell-1)!(\ell+1)!(k-\ndemi-\ell)$. Then the following identity holds 
\[L_{k}^{{\rm BG},\ell}=c_{k}^\ell L_{k}^\ell.\] 
In the critical case $\ell=\ndemi-k$, if $G_k$ is the Branson-Gover operator of \eqref{Gkbrangov} we have
\[L_k^{\rm BG}=c_kL_k, \quad G^{\rm BG}_k=(-1)^{k}c_kG_k \] 
with $c_k:=(-1)^{\ndemi-k-1} 2^{n-2k+1}((\ndemi-k)!)^2(\ndemi-k+1)=c_k^{\ndemi-k}$.\\
\textbf{(ii)} Let $Q_k$ be the operator of Definition \ref{definQ}, then 
\[Q_k^{\rm BG}=(2k-n-2)c_{k+1}Q_k.\]
\end{theo}
\textsl{Proof}: (i) For $\omega_0\in\Lambda^k(M)$, we consider the form $\omega_{F_1}$ of Lemma \ref{casnoncrit} 
of the previous section and we extend it homogeneously in a smooth $k$-form of degree $k-\ndemi+\ell$ by 
\[\begin{gathered}\til{\omega}_F=u^{k-\frac{n}{2}+\ell}\omega_{F_1}=u^{k-\frac{n}{2}+\ell}x^{\frac{n}{2}-k-\ell}\sum_{i=0}^{\ell}x^{2i}\bigl(\omega^{(t)}_i+x^2\omega_i^{(n)}\wedge\frac{dx}{x}\bigr)\\
=t^{k-\frac{n}{2}+\ell}\sum_{i=0}^{\ell}(-Q)^{i}t^{-2i}\bigl(\omega^{(t)}_i+\omega_i^{(n)}\wedge d\rho\bigr).\end{gathered}\] 

In the coordinates $u,x,y$ representing a neighbourhood $\{-1\leq Q<0, \rho<\eps\}$ and in the $k$-form bundle decomposition $\Lambda^k(\mc{C})\oplus \Lambda^{k-1}(\mc{C})\wedge\frac{du}{u}$, the exterior derivative, its dual and the form Laplacian of $\til{h}$ are given by
\begin{equation}\label{formule2}
\til{d}=\begin{pmatrix}
d & 0 \\
(-1)^ku\pl_u & d 
\end{pmatrix}, \quad \til{\delta}=u^{-2}
\begin{pmatrix}
\delta_g & (-1)^{k+1}(n+2-2k+u\pl_u) \\
0 &\delta_g 
\end{pmatrix}\end{equation}
and 
\begin{equation}\label{tildelta}
\til{\Delta}=u^{-2}\begin{pmatrix}
(u\pl_u)(u\pl_u+n-2k)+\Delta_k & 2(-1)^{k+1}d \\
2(-1)^k\delta_g & (u\pl_u-2)(u\pl_u+n-2k+2)+ \Delta_{k-1}\end{pmatrix}.
\end{equation}
So, using the properties of $\omega_{F_1}$ in Lemma \ref{casnoncrit} 
and Lemma \ref{deltaom}, we have (where $s=k-\frac{n}{2}+\ell$)
\begin{equation*}
\begin{split}
\til{\Delta}\til{\omega}_F=&\, u^{s-2}\bigl(\Delta_k+s(s+n-2k)\bigr)\omega_{F_1}+2(-1)^ku^{s-3}\delta_g\omega_{F_1}\wedge du\\
=&\,2\ell u^{s-2}x^{\ell-k+\frac{n}{2}}\bigl(L_k^\ell\omega_0+O_t(x^2)\bigr)+u^{s-2}x^{\ell-k+\frac{n}{2}+2}
\bigl(B\wedge\frac{dx}{x}+O_n(x^2)\bigr)\\
&+2(-1)^ku^{s-3}x^{\ell-k+\frac{n}{2}+2}\bigl(C+O(x^2)\bigr)\wedge du\\
=&\,(-Q)^{\ell-1}t^{k-\ell-\frac{n}{2}}\bigl(2\ell L_k^\ell\omega_0+(B+2(-1)^kC)\wedge d\rho\Bigr)+O(Q^\ell)
\end{split}
\end{equation*}
for some $(k-1)$-forms $B,C$ on $M$.
We can now apply $\ell-1$ times Lemma \ref{lem1} and get
$$
\begin{gathered}
  \til{\Delta}^\ell\til{\omega}_F=\til{\Delta}^{\ell-1}\til{\Delta}\omega_F=(-4)^{\ell-1}[(\ell-1)!]^2t^{k-\ell-\frac{n}{2}}\bigl(2\ell L_k^\ell\omega_0+(B+2(-1)^kC)\wedge d\rho\Bigr)+O(Q).
\end{gathered}
$$
Since $(n+2\til{\nabla}_T-2)$ acts on homogeneous $k$-forms of degree $k{-}\ell{-}\frac{n}{2}$ by multiplication by $-2(\ell+1)$ and
$i_T\til{d}=\mc{L}_T$ on $\mc{G}_k[\ndemi-k-\ell]$, we get
$$\mathbf{L}_k^\ell=(\ell+1)(k-\ell-\frac{n}{2})(-4)^{\ell}[(\ell-1)!]^2t^{k-\ell-\frac{n}{2}}\bigl(\ell L_k^\ell\omega_0+(\frac{B}{2}+(-1)^kC)\wedge d\rho\Bigr)$$
Note that by definition of $B_k,C_k,G_k$ we have, in the case $\ell=\frac{n}{2}-k$, \[\frac{B}{2}+(-1)^kC=(-1)^{k-1}\bigl(\frac{B_k}{2}-C_k\bigr)\omega_0=\ell G_k\omega_0.\]
(ii) Similarly, for $\omega_0\in\Lambda^{k}(M)$ closed, we set $\til{\omega}_F:=\omega'_{F_1}\wedge\demi\til{d}Q$ in $\{Q<0,\rho\leq \eps\}$ where the form $\omega'_{F_1}$ is the $0$-homogeneous expansion of $\omega'_{F_1}\in\Lambda^{k+1}(\mc{C})$ given by \eqref{F2F1}. Since $\frac{\til{d}Q}{2}=-t^2d\rho+Q\frac{dt}{t}$, we have
\begin{equation*}
\begin{split}
  \til{\omega}_F=&\,\sum_{j=0}^{\frac{n}{2}-k}x^{2j}\bigl(\omega_{2j}^{(n)}\wedge\frac{dx}{x}+x^2\omega_{2j}^{(t)}\bigr)\wedge\frac{\til{d}Q}{2}\\
=&\,\sum_{j=0}^{\frac{n}{2}-k}-(-Q)^{j}t^{-2(j-1)}\omega_{2j}^{(n)}\wedge d\rho\wedge\frac{dt}{t}+(-Q)^{j+1}t^{-2j-2}\omega_{2j}^{(t)}\wedge(-t^2d\rho+Q\frac{dt}{t})
\end{split}
\end{equation*}
and so $\til{\omega}_F$ is a smooth $(k+2)$ form. By (\ref{tildelta}) and the definition of $B'_k,D'_k$ we have
\begin{equation*}
\begin{split}
\til{\Delta}\til{\omega}_F=&\,\til{\Delta}(-u^2\omega_{F_1}'\wedge\frac{du}{u})=2(-1)^kd\omega'_{F_1}-(\Delta_{k+1}\omega'_{F_1})\wedge\frac{du}{u}\\
=&\,-x^{n-2k-2}2D'_k\omega_0\wedge\frac{dx}{x}+(-1)^{k+1}x^{n-2k}\bigl(B'_k\omega_0\wedge\frac{dx}{x}+\omega_1\bigr)\wedge\frac{du}{u}+O(Q^{\ndemi-k})\\
=&\,(-1)^{\frac{n}{2}-k-1}2Q^{\frac{n}{2}-k-2}t^{2k-n+4}D'_k\omega_0\wedge d\rho\\
&+(-1)^{\frac{n}{2}+1}Q^{\frac{n}{2}-k-1}t^{2k-n+2}\bigl(B'_k\omega_0\wedge\frac{dt}{t}+(-1)^k\omega_1\bigr)\wedge d\rho+O(Q^{\frac{n}{2}-k})
\end{split}\end{equation*}
for some form $\omega_1$ on $M$, the value of which is not important for our purpose.
By Lemma \ref{lem1}, we have
\begin{equation*}
\begin{split}
  \til{\Delta}^2\til{\omega}_F=&\,(-1)^{\frac{n}{2}-k-1}2Q^{\frac{n}{2}-k-2}t^{2k-n+4}\til{\Delta}(D'_k\omega_0\wedge d\rho)\\
&+4(\frac{n}{2}-k-1)(-1)^{\frac{n}{2}+1}Q^{\frac{n}{2}-k-2}t^{2k-n+2}\bigl(B'_k\omega_0\wedge\frac{dt}{t}+(-1)^k\omega_1\bigr)\wedge d\rho+O(Q^{\frac{n}{2}-k-1})
\end{split}\end{equation*}
and by \eqref{tildelta}, we have
\begin{equation*}
\begin{split}
  \til{\Delta}(D'_k\omega_0\wedge d\rho)=&\,\til{\Delta}(x^2D'_k\omega_0\wedge\frac{dx}{x})\\
=&\,u^{-2}\Delta_{k+2}(x^2D'_k\omega_0\wedge\frac{dx}{x})+2u^{-2}(-1)^k\delta_g(x^2D'_k\omega_0\wedge\frac{dx}{x})\wedge\frac{du}{u}\\
=&\,2(-1)^kt^{-2}\delta_0 D'_k\omega_0\wedge d\rho\wedge\frac{dt}{t}+2(2k-n-4)t^{-2}D'_k\omega_0\wedge\frac{dt}{t}
\end{split}\end{equation*}
where we have used \eqref{formule}, \eqref{decompP} and $dD'_k\omega_0=0$. We thus have
\begin{equation*}
\begin{split}
\til{\Delta}^2\til{\omega}_F=&\,Q^{\frac{n}{2}-k-2}t^{2k-n+2}\Bigl[4(-1)^{\frac{n}{2}-1}\Bigl(\delta_0 D'_k\omega_0-(\frac{n}{2}-k-1)B'_k\omega_0\Bigr)\wedge d\rho\wedge\frac{dt}{t}+\omega'_1\wedge\frac{dt}{t}\\
&\,\quad\quad\quad\quad\quad\quad\quad\quad+\omega'_2\wedge d\rho\Bigr]+O(Q^{\frac{n}{2}-k-1})\\
=&\,Q^{\frac{n}{2}-k-2}t^{2k-n+2}\Bigl[4(n-2k)(\frac{n}{2}-k-1)(-1)^{\frac{n}{2}+k}Q_k\omega_0\wedge d\rho\wedge\frac{dt}{t}+\omega'_1\wedge\frac{dt}{t}\\
&\,\quad\quad\quad\quad\quad\quad\quad\quad+\omega'_2\wedge d\rho\Bigr]+O(Q^{\frac{n}{2}-k-1})\\
\end{split}\end{equation*}
by Corollary \ref{cor1}, and where $\omega'_1,\omega'_2$ are forms in $\Lambda^{k+1}(M)$. By iterative use of Lemma \ref{lem1}, we get
\begin{equation*}
\begin{split}
\til{\Delta}^{\frac{n}{2}-k}\til{\omega}_F=&\,\til{\Delta}^{\frac{n}{2}-k-2}\til{\Delta}^2\til{\omega}_F\\
=&\,t^{2k-n+2}\Bigl[2^{n-2k-1}(\frac{n}{2}-k)[(\frac{n}{2}-k-1)!]^2(-1)^{\frac{n}{2}+k}Q_k\omega_0\wedge d\rho\wedge\frac{dt}{t}\\
&\quad\quad\quad\quad\quad+\omega'_1\wedge\frac{dt}{t}+\omega_2'\wedge d\rho\Bigr]+O(Q)
\end{split}\end{equation*}
we infer from the definition of $Q_{k}^{\rm BG}$ that
$$Q_k^{BG}=(-1)^{\frac{n}{2}+k+1}2^{n-2k}(\frac{n}{2}-k+1)!(\frac{n}{2}-k-1)!Q_k$$
\qed

\section{Proof of the main results}

We start with the proof of Theorem \ref{theo2}.

\textsl{Proof of Theorem \ref{theo2}}: the existence of $\omega$ in (i) is proved in Proposition \ref{poisson}.
The fact that the log terms $L_k,Q_k$ coincide with the Branson-Gover operators follows from Theorem \ref{coincide}.
The uniqueness of the solution is rather clear by construction: using the arguments used in the proof of Proposition
\ref{poisson}, a solution in $C^{\ndemi-k,\alpha}(\bar{X},\Lambda^k(\bar{X}))$ would have its 
first $\ndemi-k$ Taylor coefficients uniquely (and locally) determined 
by the boundary value $\omega_0$ and then two such solutions with same boundary data 
would agree to order $x^{\ndemi-k+\alpha}$ and would then be in 
$L^2(X,\Lambda^k(X))$. The proof of (ii) is similar and follows from Proposition \ref{construQ} and Theorem \ref{coincide}.
\qed\\

\textsl{Proof of Theorem \ref{th1}}: The infinite dimensionality of $K_m^k(\bar{X})$ for $m<n-2k+1$
follows from Proposition \ref{poisson}. Indeed for $m<n-2k$ this is clear since the solution of \eqref{dirichletpb}
are parameterized by $C^{\infty}(M,\Lambda^{k}(M))$. If $m=n-2k$,
one can use that there is an infinite set of $\omega_0\in C^{\infty}(M,\Lambda^{k}(M))$ 
such that $G_k\omega_0\not=0$ and $L_k\omega_0=0$ since $\ker L_k$ is infinite dimensional and 
$\ker G_k\cap\ker L_k$ is finite dimensional by ellipticity of $dG_k+L_k$. Solutions of \eqref{dirichletpb} are then in $C^{n-2k}(\bar{X},\Lambda^k(\bar{X}))$.

The finite dimensionality for $m= n-2k+1$ is a little more involved. 
Let $\omega$ be a harmonic form in $C^{n-2k+1}(\bar{X},\Lambda^{k}(\bar{X}))$, then
Taylor expanding, there exist some forms $\omega_j^{(n)},\omega_j^{(t)}\in C^{n-2k+1-j}(M,\Lambda(M))$ so that
\[\omega-\sum_{j=0}^{n-2k}x^j(\omega_j^{(t)}+\omega_j^{(n)}\wedge dx)\in x^{n-2k+1}L^\infty(\Lambda^k(\bar{X})),\]
and $L_k\omega_0=0$. Now by Lemma \ref{maz2} we know that $\omega$ has a weak expansion to order $x^N$
with values in $H^{-N}(M)$ like in \eqref{weakexp} for any $N>0$ large. 
Moreover $\delta_g\omega$ is also a harmonic form in $C^{n-2k}(\bar{X},\Lambda^{k-1}(M))$
which is a $O(x)$ and has an expansion to order $x^{N}$ with values in $H^{-N-1}(M)$ for any $N$. 
Now, using the indicial equation like in the proof of Proposition \ref{poisson}, 
the weak expansion of $\delta_g\omega$ vanish to order $x^{n-2k+2}$, so in particular we obtain
$\delta_g\omega\in x^{n-2k}L^\infty(\Lambda^{k-1}(\bar{X}))$ from the regularity of $\omega$. 
Then $\delta_g\omega\in 
L^2(\Lambda^{k-1}(X))$ for $k<\ndemi-1$, while for $k=\ndemi-1$ it is in $L^2$ if we assume in 
addition that $\omega\in C^{n-2k+1,\alpha}(\bar{X},\Lambda^k(\bar{X}))$ for some $\alpha>0$ (since
then $\delta_g\omega\in x^{n-2k+\alpha}L^\infty(\Lambda^{k}(\bar{X}))$). 
But as shown in the proof of Proposition \ref{poisson}, an $L^2$ harmonic form which is 
coclosed is identically $0$. Now we can apply the result of Proposition \ref{formalequations} (see the Remark below Corollary \ref{cor2}), and compute $\delta_g\omega$, which gives $G_k\omega_0=0$. Since $dG_k+L_k$ is elliptic, $\ker L_k\cap\ker G_k$ is finite dimensional and contains only smooth forms, so $\omega_0$ is smooth. Then $\omega$ is polyhomogeneous and is the solution of Proposition \ref{poisson}, up to 
an element of $\ker_{L^2}(\Delta_k)$, it is then in $C^{n-1}(\bar{X},\Lambda^k(\bar{X}))$ in general and in 
$C^{\infty}(\bar{X},\Lambda^{k}(\bar{X}))$ if $(X,g)$ smooth Poincar\'e-Einstein manifold.

Let $m\in[n-2k+1,n-1]$ be an integer. The exact sequence \eqref{shortexact} is defined by inclusion of $\iota: H^{k}(\bar{X},\pl\bar{X})\to K_{m}^k(\bar{X})$
and restriction to the boundary $r:K_m^k(\bar{X})\to \mc{H}_L^k(M)$, here of course we use
the identification $H^{k}(\bar{X},\pl\bar{X})\simeq \ker_{L^2}(\Delta_k)$ and the regularity of harmonic $L^2$ 
forms in Theorem \ref{Mazzeo}. 
The injectivity of $\iota$ is clear, the surjectivity of $r$ comes from Proposition
\ref{poisson}, the definition of $\mc{H}_L^k$ and Theorem \ref{coincide}. The kernel of $r$ is composed of those forms
of $K_m^k(\bar{X})$ which vanish at $M$, but by Proposition \ref{poisson}, these are $L^2$, and thus in the image of 
$H^{k}(\bar{X},\pl\bar{X})$ by the map $\iota$.
\qed\\  

\textsl{Proof of Theorem \ref{th3}}: First note that the space $Z^k(\bar{X})$ in Theorem \ref{th3} is included in $K_{n-2k+1}^k(\bar{X})$, and thus  of finite dimension and composed of forms in $C^{n-1}(\bar{X},\Lambda^{k}(\bar{X}))$ (even in the case $k=\ndemi$ by the arguments above).

(i) the maps in the complex
\[ 0\xrightarrow{}H^{k}(\bar{X},\pl\bar{X})\xrightarrow{\iota} Z^k(\bar{X})\xrightarrow{r}
\mc{H}^k(\pl\bar{X})\xrightarrow{d_e}H^{k+1}(\bar{X},\pl\bar{X})  \]
are defined  as follows: $\iota$ is given by 
inclusion where $H^{k}(\bar{X},\pl\bar{X})\simeq \ker_{L^2}(\Delta_k)$, this is well defined since
$L^2$ harmonic forms are closed, coclosed and in $C^{n-2k+1}(\bar{X},\Lambda^{k}(\bar{X}))$; 
$r$ is defined as restriction at the boundary and it maps in $\mc{H}^k(M)$ since  $r(\omega)\in\ker L_k\cap\ker G_k$ by the discussion above
and $d\omega=0$ implies $dr(\omega)=0$; 
the last map $d_e$ is the composition $d_e=d\circ \Phi$ where $\Phi: C^{\infty}(M,\Lambda^k(M))\to C^{\infty}(X,\Lambda^k(X))/\ker_{L^2}(\Delta_k)$
is defined by $\Phi(\omega_0)=\omega$ where $\omega$ is the solution of \eqref{dirichletpb} in Proposition 
\ref{poisson}. Note that $\Phi$ is only defined modulo $\ker_{L^2}(\Delta_k)$ and linear by uniqueness of the solution
in \eqref{dirichletpb} modulo $\ker_{L^2}(\Delta_k)$. Applying $d$ kills the indeterminacy with respect to $\ker_{L^2}(\Delta_k)$ since $L^2$ harmonic forms are closed. Then $d\Phi(\omega_0)$ is harmonic and
since the boundary value of $\Phi(\omega_0)$ is closed, then $d\Phi(\omega_0)=O(x)$, and by Proposition 
\ref{poisson} it is in $L^2$.
For the exactness of the sequence, first note that $\ker r$ is composed of closed and coclosed forms  
which are $O(x)$, this implies that those forms are $L^2$ by Proposition \ref{poisson}, so $\textrm{Im}\, \iota=\ker r$
since also $L^2$ harmonic forms vanish at the boundary. Now $\omega_0\in \ker d_e$ if 
$\Phi(\omega_0)$ is closed, but it is also coclosed and in 
$C^{n-2k+1}(\bar{X},\Lambda^{k}(\bar{X}))$ by Proposition \ref{poisson} and the fact that $\omega_0\in \ker d\cap\ker G_k\subset \ker L_k\cap\ker G_k$,
therefore $\Phi(\omega_0)\in Z^k(\bar{X})$ and $\omega_0\in\textrm{Im}\,r$. Moreover by Proposition \ref{poisson} we have $\Phi(r(\omega))-\omega\in\ker_{L^2}(\Delta_k)$, this implies $\textrm{Im}\,r\subset \ker d_e$, this proves exactness of the sequence.

(ii) the map in the complex \eqref{longcomplex} are defined similarly:
first $\iota:H^{k}(\bar{X},\pl\bar{X})\to [Z^k(\bar{X})]$ is the composition of the inclusion
$\ker_{L^2}(\Delta_k)\to Z^k(\bar{X})$ with the natural map $Z^k(\bar{X})\to [Z^k(\bar{X})]$ obtained
by taking cohomology class. The map $r:[Z^k(\bar{X})]\to [\mc{H}^k(\pl\bar{X})]$ is the map induced by 
the restriction map $Z^k(\bar{X})\to \mc{H}^k(\pl\bar{X})$ used in (i). This is well defined since
if $d\alpha\in Z^k(\bar{X})$, then $r(d\alpha)=d\alpha_0$ where $\alpha_0=\alpha|_{\pl\bar{X}}$, and so
$[r(d\alpha)]=0$ if $[\,\cdot\,]$ denotes cohomology class in $H^k(\pl\bar{X})$.
The last map $d_e:[Z^k(\pl\bar{X})]\to H^{k+1}(\bar{X},\pl\bar{X})$ is the map induced by $d_e$ defined in (i), i.e.
$d_e=d\circ\Phi$ where $\Phi$ maps $\omega_0$ to the solution of \eqref{dirichletpb}. Note that it is well defined since
for $d\alpha_0\in \mc{H}^k(\pl\bar{X})$, we have $d_e(d\alpha_0)=d\Phi(d\alpha_0)$ and, by uniqueness 
of the solution of \eqref{dirichletpb}, $\Phi(d\alpha_0)-d\Phi(\alpha_0)\in\ker_{L^2}(\Delta_{k+1})$ thus
$d\Phi(d\alpha_0)=0$. 

To show that $\ker r={\rm Im}\, \iota$, we need to show that if $\omega \in Z^k(\bar{X})$ is a representative in 
$[Z^k(\bar{X})]$ such that $r(\omega)=d\alpha_0$ for some smooth $\alpha_0$, then there is $\omega'\in\ker_{L^2}(\Delta_k)$
such that $\omega-\omega'$ is exact.
But as said above, we have $\Phi(d\alpha_0)-d\Phi(\alpha_0)\in\ker_{L^2}(\Delta_{k})$ and $\Phi(r(\omega))-\omega\in\ker_{L^2}(\Delta_k)$ thus $\omega-d\Phi(\alpha_0)\in \ker_{L^2}(\Delta_k)$ and we are done. 
To show that $\ker d_e={\rm Im}r$, we need to prove that for $\omega_0\in \mc{H}^k(\pl\bar{X})$ 
a representative in $[\mc{H}^k(\pl\bar{X})]$ then $\Phi(\omega_0)$ is closed if and only if there exists
$\omega\in Z^k(\bar{X})$ so that $r(\omega)-\omega_0$ is exact. But $\Phi(\omega_0)$ is in $Z^k(\bar{X})$ if $d\Phi(\omega_0)=0$, 
thus $\ker d_e\subset {\rm Im}\, r$; conversely if there is $\omega\in Z^k(\bar{X})$ with $\omega=\omega_0+d\alpha_0+O(x)$,
then $\omega-\Phi(\omega_0+d\alpha_0)\in \ker_{L^2}(\Delta_k)$ and so $d\Phi(\omega_0)=0$ since $\Phi(d\alpha_0)-d\Phi(\alpha_0)\in\ker_{L^2}(\Delta_{k})$. To conclude, we need to prove that
${\rm Im}\, d_e\subset \ker\iota$. But this is clear since $d_e\omega_0=d\Phi(\omega_0)$ is an exact $(k+1)$-form in $L^2$
with $\Phi(\omega_0)\in C^{n-2k+1}(\bar{X},\Lambda^{k}(\bar{X}))$. Note that in the case $k=\ndemi$, we make use of Proposition \ref{nader}.
 
(iii) Suppose that $[\mc{H}^k(\pl\bar{X})]=H^k(\pl\bar{X})$. If $\omega\in\ker\iota$, it is a $k$-form in $\ker_{L^2}(\Delta_k)$
which can be written $\omega=d\alpha$ with $\alpha$ smooth.
Moreover if $\alpha_0=\alpha|_{\pl\bar{X}}$, then $d(\Phi(\alpha_0)-\alpha)\in\ker_{L^2}(\Delta_k)$ and 
$\Phi(\alpha_0)-\alpha=O(x)$, an easy integration by part shows that $d\Phi(\alpha_0)=d\alpha=\omega$.
Here $\alpha_0$ is closed since $\omega=O(x)$, but by assumption there is a $\alpha_0'\in\mc{H}^k(\pl\bar{X})$
such that $\alpha_0-\alpha_0'=d\beta$ for some smooth $\beta$. Since now $d\Phi(d\beta)=d[\Phi,d]\beta=0$,
we have $d_e\alpha_0'=\omega$ and $\omega\in {\rm Im}\,d_e$, which gives $\ker\iota={\rm Im}\,d_e$. Eventually, the equality $[Z^k(\bar{X})]=H^k(\bar{X})$ is clear from the discussion above since $[Z^k(\bar{X})]\subset H^k(\bar{X})$ and  
\[\begin{gathered}
H^{k}(\bar{X},\pl\bar{X})\xrightarrow{\iota}[Z^k(\bar{X})]\xrightarrow{r}H^k(\pl\bar{X})
\xrightarrow{d_e}H^{k+1}(\bar{X},\pl\bar{X})\\
H^{k}(\bar{X},\pl\bar{X})\xrightarrow{\iota}H^k(\bar{X})\xrightarrow{r}H^k(\pl\bar{X})
\xrightarrow{d_e}H^{k+1}(\bar{X},\pl\bar{X})
\end{gathered}\]
are both exact sequences. 

As for the converse, if $\ker \iota^{k+1}={\rm Im}\,d_e^k$ and $[Z^k(\bar{X})]=H^{k}(\bar{X})$, then we have the exact sequences
\[\begin{gathered}
H^{k}(\bar{X})\xrightarrow{r}[\mc{H}^k(M)]\xrightarrow{d_e}H^{k+1}(\bar{X},M)\xrightarrow{\iota}[Z^{k+1}(\bar{X})]\\
H^{k}(\bar{X})\xrightarrow{r}H^k(M)\xrightarrow{d_e}H^{k+1}(\bar{X},M)
\xrightarrow{\iota'}H^{k+1}(\bar{X})
\end{gathered}\]
and since $[Z^{k+1}(\bar{X})]\subset H^{k+1}(\bar{X})$, we obviously have
$\ker\iota=\ker\iota'={\rm Im}\, d_e$ and so $[\mc{H}^k(M)]= H^k(M)$ (recall $[\mc{H}^k(M)]\subset H^k(M)$).
\qed\\

\textsl{Proof of Proposition \ref{qpositif}}: Assume $\cjg Q_kv,v\cjd\geq 0$. 
To show surjectivity of $\mc{H}^k(M)\to H^k(M)$, we need to prove that for all $\omega_0\in C^{\infty}(M,\Lambda^k(M))$ closed,
there exists an exact form $d\alpha$ (with $\alpha\in C^{\infty}(M,\Lambda^k(M))$) such that
$G_k(\omega_0+d\alpha)=0$. Consider $\Box:=\delta_0Q_{k}d+(d\delta_0)^{\ndemi-k+1}$ which is elliptic, self-adjoint 
and non-negative if $Q_{k}\geq 0$. Its kernel is finite dimensional (containing $\ker(d+\delta_0)$) 
and all $v\in \ker\Box$ are smooth by elliptic regularity, and satisfy $\cjg \delta_0Q_kdv,v\cjd_{L^2}=0$, which implies
$\cjg Q_kdv,dv\cjd_{L^2}=0$. 
Let $\textbf{H}\subset L^2(\Lambda^{k}(M))$ be the $L^2$ completion of the set $C^{\infty}(M,\Lambda^{k}(M))\cap \ker d$ of smooth closed forms 
and let us define $\textbf{Q}$ the symmetric form $\textbf{Q}(v,v):=
\cjg Q_kv,v\cjd_{L^2}$ on $\textbf{H}$,   
it is a non-negative form induced by $\Pi_{\textbf{H}}Q_k$ on $\textbf{H}$ where $\Pi_{\textbf{H}}$ denotes orthogonal projection
from $L^2(\Lambda^{k}(M))$ to $\textbf{H}$. The form has a domain $D(\textbf{Q})$ and Friedrichs extension theorem implies
that there exists a self adjoint operator $Q^{\rm Fr}_k:\textbf{H}\to \textbf{H}$ with domain $D(Q^{\rm Fr})$ 
such that $\cjg Q^{\rm Fr}_ku,u\cjd=\textbf{Q}(u,u)$ for $u\in D(\textbf{Q})\cap D(Q^{\rm Fr})$.
But clearly $d(C^{\infty}(M,\Lambda^{k-1}(M)))\subset D(Q^{\rm Fr}_k)$ and so $\Pi_{\textbf{H}}Q_kdv=Q^{\rm Fr}_kdv$
for $v$ smooth. Using now the spectral theorem for $Q^{\rm Fr}_k$, we see that $Q^{\rm Fr}_kdv=0$ with $v$ smooth 
if and only if $\cjg Q_kdv,dv\cjd=0$ and $v$ is smooth, thus in particular if $v\in \ker \Box$. Thus 
$Q_kdv\perp \omega$ for all $\omega\in \textbf{H}$ if $v\in\ker \Box$.
Now this implies that, with $\omega$ closed and smooth, we have $\cjg v,G_k\omega\cjd=\cjg Q_kdv,\omega\cjd=0$ 
for $v\in\ker\Box$ since $Q_k$ is symmetric on closed forms, and so $G_k\omega$ is in the range of $\Box$
and there exists $\alpha$ such that $\Box \alpha=-G_k\omega$, but since ${\rm Im}\,G_k\subset {\rm Im}\,\delta_0$
which is orthogonal to ${\rm Im}\,d$, we deduce that $(d\delta_0)^{\ndemi-k+1}\alpha=0$ and this achieves the proof.
Note in particular that in this case $\{d\varphi; L_{k-1}\varphi=0\}=\{d\varphi; Q_kd\varphi\in{\rm Im}\,\delta_0\}$,
see Corollaries 2.12 and 2.13 of \cite{BG} for discussions about these spaces.
\qed\\

\section{Computations in some special cases}
\label{Comput}

In this section we compute the operator $L_k$, $G_k$ and $Q_k$ in dimension $4$ and $6$.

\begin{lem}\label{comp}
  Let $(M^4,h)$ a four dimensional Riemannian manifold and define for any symmetric $2$-tensor $H$ the map 
  $j(H):=J(h^{-1}H)$ where $J$ is defined in \eqref{jh}. Then we have
$$\displaylines{
\hfill L_1=\frac{1}{2}\delta d,\hfill G_1=-\frac{1}{4}\delta\Bigl(\Delta-2j(\Ric)+\frac{2}{3}\Scal\Bigr),\hfill Q_1=\frac{1}{2}\Bigl(\Delta-2j(\Ric)+\frac{2}{3}\Scal\Bigr),\cr
\hfill L_0=-\frac{1}{16}\delta\Bigl(\Delta-2j(\Ric)+\frac{2}{3}\Scal\Bigr)d,\hfill G_0=0,\hfill Q_0=-\frac{1}{24}\Bigl(\Delta\Scal-3|\Ric|^2+\Scal^2\Bigr)}$$
where $\rm{Ric}$ is the Ricci tensor of $h$ and $\rm{Scal}$ its scalar curvature
\end{lem}

\textbf{Remark}: If $n=4$, $L_{\ndemi-2}$ is the Paneitz operator (up to a constant factor). The result of Gursky and 
Viaclovsky \cite{Gur}
says that if the Yamabe invariant $Y(M,[h_0])$ is positive and  
\[\int_M Q{\rm dvol}_{h_0}+\frac{1}{6}Y(M,[h_0])^2>0\]
then $L_0$ is a non-negative operator with kernel reduced to constants. Combining with Theorem 2.6 of Branson-Gover\cite{BG}, we have that $\mc{H}^1(M)\simeq H^1(M)$ and there is a conformally invariant basis of $H^1(M)$ with respect to $[h_0]$
made of conformal harmonics.

Using the inequality $\|D\omega\|_2^2\geq\|\delta\omega\|^2/n$ for all $1-$form $\omega$ and the Bochner formula we get
\begin{cor}
  Let $M^4$ be a four dimensional manifold and $\lambda_1(x)\geq\cdots\geq\lambda_4(x)$ the eigenvalues of its Ricci curvature at $x$. If $\lambda_2(x)+\lambda_3(x)+\lambda_4(x)\geq0$ for all $x\in M$ then $\mc{H}^1(M)\to H^{1}(M)$ is surjective.
\end{cor}

\textsl{Proof}: 
  For any closed form $\omega$, we have $\langle\Delta\omega,\omega\rangle=\|\delta\omega\|_2^2=\|D\alpha\|_2^2+\Ric(\omega,\omega)\geq\|\delta\omega\|^2/4+\int_M\Ric(\omega,\omega)$, and so $\langle\Delta\omega,\omega\rangle\geq\frac{4}{3}\int_M\Ric(\omega,\omega)$.
$$\langle Q_1\omega,\omega\rangle\geq\frac{1}{3}\int_M\Scal|\omega|^2-\Ric(\omega,\omega)$$
\qed

\begin{lem}\label{comp1}
  Let $(M^6,h)$ a six dimensional manifold. If $j$ is defined like in Lemma \ref{comp} and $\tr(H)$ denotes the trace of $H$ with respect to $h$, 
  then we have
$$\displaylines{
\hfill L_2=\frac{1}{2}\delta d,\hfill G_2=\frac{1}{4}\delta\Bigl(\Delta-j(\Ric)+\frac{2}{5}\Scal\Bigr),\hfill Q_2=\frac{1}{2}\Bigl(\Delta-j(\Ric)+\frac{2}{5}\Scal\Bigr),\cr
 L_1=-\frac{1}{16}\delta\Bigl(\Delta-j(\Ric)+\frac{2}{5}\Scal\Bigr)d,\cr
 G_1=\frac{1}{16}\Bigl[\delta\Delta^2-\frac{\delta d\delta}{2}j(\Ric-\frac{3}{10}\Scal)-\delta j(2\Ric-\frac{3}{5}\Scal)\Delta-\frac{\delta d}{20}\Scal\delta\hfill\cr
\hfill+\delta j\bigl(2B-\tr(B)+\frac{3\Ric^2}{4}-\frac{16\Scal \,\Ric}{5}+\frac{449\Scal^2}{100}\bigr)\Bigr],\cr
Q_1=-\frac{1}{4}\Bigl[\Delta^2-\frac{d\delta}{2}j(\Ric-\frac{3}{10}\Scal)-j(2\Ric-\frac{3}{5}\Scal)\Delta-\frac{ d\Scal\delta}{20}\hfill\cr
\hfill+ j\bigl(2B-\tr(B)+\frac{3\Ric^2}{4}-\frac{16\Scal \,\Ric}{5}+\frac{449\Scal^2}{100}\bigr)\Bigr],\cr
L_0=\frac{1}{96}\Bigl[(\delta d)^3-\frac{\delta d\delta}{2}j(\Ric-\frac{3}{10}\Scal)d-\delta j(2\Ric-\frac{3}{5}\Scal)d\delta d-\frac{\delta d}{20}\Scal\delta d\hfill\cr
\hfill+\delta j\bigl(2B-\tr(B)+\frac{3\Ric^2}{4}-\frac{16\Scal \,\Ric}{5}+\frac{449\Scal^2}{100}\bigr)d\Bigr],\cr
 G_0=0,\cr
Q_0=\frac{1}{640}\Bigl[\Delta^2\Scal+\Scal\Delta\,\Scal+2(\Ric,{\rm Hess}\,\Scal)-20\Delta\tr(B)-40\Delta|P|^2\hfill\cr
\hfill+\frac{2}{25}\Scal^3-12\,\Scal\,\tr(B)-80\,\tr (P^3)-80(P,B)\Bigr],}$$
where $B$ denotes the Bach tensor of $h$, $P$ its Schouten tensor, $\rm{Ric}$ its Ricci tensor and $\rm{Scal}$ its scalar curvature.
\end{lem}

\begin{lem}\label{Hn/2-1}
For any $n\geq4$, we have the identities\begin{eqnarray*}
G_{\frac{n}{2}-1}&=&(-1)^{\frac{n}{2}+1}\Bigl(\frac{\delta d\delta}{4}-\frac{\delta j(P)}{2}+\delta\frac{\tra(P){\rm Id}}{4}\Bigr),\\
&=&(-1)^{\frac{n}{2}+1}\frac{\delta d\delta}{4}+(-1)^\frac{n}{2}\delta\Bigl(\frac{j(\Ric)}{n-2}-\frac{\Scal{\rm Id}}{2(n-1)}\Bigr)\\
L_{\frac{n}{2}-2}&=&-\delta\Bigl(\frac{d\delta}{16}-\frac{j(\Ric)}{4(n-2)}+\frac{\Scal\,{\rm Id}}{8(n-1)}\Bigr)d\\
Q_{\frac{n}{2}-1}&=&\Bigl(\frac{\Delta}{2}-\frac{2j(\Ric)}{n-2}+\frac{\Scal\,{\rm Id}}{n-1}\Bigr)
\end{eqnarray*}
\end{lem}

For the non critical case, we have
\begin{lem}\label{comp2}
We set $j^{\sharp}(H)=2j(H)-\tr(H)\,{\rm Id}$. For any $n\geq3$, we have
\begin{eqnarray*}
L^1_{k}&=&\frac{\delta d}{2}+\frac{(n-2k-2)d\delta}{2(n-2k+2)}+\frac{(n+k-2)(n-2k-2)}{8(n-1)(n-2)}\Scal-\frac{(n-2k-2)j(\Ric)}{2(n-2)}
\end{eqnarray*}
which generalizes the conformal Laplacian on functions,
\begin{eqnarray*}
L^2_{k}&=&-\frac{n-2k-4}{16}\Bigl(\frac{(d\delta)^2}{n-2k+4}+\frac{(\delta d)^2}{n-2k-4}+\frac{2d j^{\sharp}(P)\delta}{n-2k+4}-\frac{2\delta j^{\sharp}(P)d}{n-2k-4}\\
& &-\frac{j(P)\Delta+\Delta j^{\sharp}(P)}{2}+j^{\sharp}(P^2+\frac{B}{n-4})+\frac{(n-2k)j^{\sharp}(P)^2}{4}\Bigr)
\end{eqnarray*}
which generalizes the Paneitz-Branson operator on functions.
\end{lem}
\noindent\textsl{Proofs of Lemmas \ref{comp}, \ref{comp1}, \ref{Hn/2-1} and \ref{comp2}}: This is a quite tedious computation, therefore we do not give the full details. By \cite[Eq. (3.18)]{FG2}, we have
$$h_x=h_0-x^2P+x^4\frac{h_2}{8}-x^6\frac{h_3}{48}+o(x^6),$$
where $P=\frac{1}{n-2}\bigl(\Ric-\frac{\Scal}{2(n-1)}\bigr)$, $h_2=-\frac{2B}{n-4}+2P^2$ and $\tr(h_3)=-\frac{8\tr(PB)}{n-4}$ 
and in the case $n=4$ we take $B=h_3=0$; note that we have ignored the first log term in the metric expansion (i.e. the obstruction tensor) 
in dimension $4$ and $6$ since, as it is clear from Lemma \ref{formlapl},
they do not show up in the construction the $L_k^\ell,G_k,Q_k$. We set $B'=\frac{2B}{n-4}$, then, with the notations of the proof of Lemma \ref{formlapl}, we have 
$$
\displaylines{
O_x=I+x^2\frac{P}{2}+x^4\frac{4P^2+B'}{16}+x^6\frac{h_3+12P^3+5PB'+4B'P}{96}+o(x^6),\cr
I_x=I+x^2\frac{J(P)}{2}+x^4\frac{J(2P^2+B')+2J(P)^2}{16}\hfill\cr
\hfill+x^6\bigl(\frac{J(h_3+4P^3+5PB'+B'P)+3J(P)J(B')+6J(P^2)J(P)+2J(P)^3}{96}\bigr)+o(x^6).}$$
Then 
$$\displaylines{\star_x=\star_0+x^2\frac{[\star_0,J(P)]}{2}+x^4\Bigl(\frac{2\bigl[J(P),[J(P),\star_0]\bigr]-[J(2P^2+B'),\star_0]}{16}\Bigr)\hfill\cr
\hfill+ x^6\Bigl(\frac{\bigl[4J(P)^3-J(h_3+4P^3+5PB'+B'P),\star_0\bigr]+3\bigl[J(B'),[J(P),\star_0]\bigr]}{96}\Bigr)}$$
from which we infer that
$$\displaylines{\star_x^{-1}\bigl[\pl_x,\star_x\bigr]=x\star_0^{-1}[\star_0,J(P)]+x^3\frac{\star_0^{-1}[\star_0,J(2P^2+B')]}{4}\hfill\cr
\hfill+ x^5\Bigl(\frac{\star_0^{-1}\bigl[4J(P)^3-J(h_3+4P^3+6PB'),\star_0\bigr]-2(\star_0^{-1}[\star_0,J(P)])^3-6\star_0^{-1}[\star_0,J(P)]\star_0^{-1}[\star_0,J(P^2)]}{16}\Bigr),\hfill\cr
\delta_x=\delta_0+x^2\frac{\bigl[\delta_0,\star_0^{-1}[\star_0,J(P)]\bigr]}{2}+x^4\frac{\bigl[\delta_0,\star_0^{-1}[\star_0,J(2P^2+B')]\bigr]+2\Bigl[\star_0^{-1}[\star_0,J(P)],\bigl[\star_0^{-1}[\star_0,J(P)],\delta_0\bigr]\Bigl]}{16},\cr
\Delta_k=
\begin{pmatrix}
  -(x\pl_x)^2+(n-2k)x\pl_x&2(-1)^{k+1}d\\
0&-(x\pl_x)^2+(n-2k+2)x\pl_x
\end{pmatrix}\hfill\cr
\hfill+x^2
\begin{pmatrix}
  \Delta_0-(2J(P)-\tr P)x\pl_x&(-1)^k[d,2J(P)-\tr P]\\
2(-1)^{k+1}\delta_0&\Delta_0-(2J(P)-\tr P)(2+x\pl_x)
\end{pmatrix}\cr
+x^4
\begin{pmatrix}
  A_1-A_2x\pl_x&(-1)^k[d,A_2]\\
(-1)^k2[2J(P)-\tr P,\delta_0]&A_1-A_2(4+x\pl_x)
\end{pmatrix}\cr
+x^6 
\begin{pmatrix}
  A_3&A_4\\
A_5&A_6
\end{pmatrix}+o(x^6)}$$
where $A_1=\frac{d[\delta_0,2J(P)-\tr P]+[\delta_0,2J(P)-\tr P]d}{2}$, $A_2=J(P^2+\frac{B'}{2})-\frac{1}{2}\tr(P^2+\frac{B'}{2})$ and $A_61=\frac{9}{40}|P|^2\Scal-\frac{3}{2}\tr (P^3)-\frac{3}{4}g(P,B')$.

For $n=6$ and $k=1$ we follow the formal method of Subection \ref{theoperatorqk} and find
$$\omega_{F_2}'=\frac{dx}{x}-x^2(\frac{d\Scal}{80}+\frac{\Scal dx}{40x})+x^4(\frac{\Delta\Scal}{160}+\frac{\Scal^2}{800}-\frac{\tr B}{8}-\frac{|P|^2}{4})\frac{dx}{x}$$
and so by computing $\Delta_k\omega_{F_2}$ one finds
$$\displaylines{Q_0=\frac{1}{640}\Bigl(\Delta^2\Scal+\Scal\Delta\Scal+2(\Ric,{\rm Hess}\,\Scal)-20\Delta\tr(B)-40\Delta|P|^2\hfill\cr
\hfill+\frac{2}{25}\Scal^3-12\Scal\,\tr(B)-80\tr (P^3)-80(P,B)\Bigr)}$$
The other computations are made by the same way. For instance, for $k=n/2-1$, we have
\[\Delta \omega_{F_1}= x^2\delta_0 d\omega_0+x^3(-1)^{\frac{n}{2}+1}\bigl(\frac{\delta_0d\delta_0\omega_0}{2} -2\delta_0A\omega_0\bigr)\wedge dx+O(x^4),\]
and so $$B_{\frac{n}{2}-1}\omega_0=-\frac{\delta_0d\delta_0\omega_0}{2}+2\delta_0A\omega_0$$
We have $\delta\omega_{F_1}=\frac{x^4}{2}\delta_0A\omega_0+O(x^5)$, and so $C_{\frac{n}{2}-1}=\frac{\delta_0 A}{2}$.
By Proposition \ref{hklocal}, we have that
\[G_{\frac{n}{2}-1}=(-1)^{\frac{n}{2}+1}\Bigl(\frac{\delta_0d\delta_0}{4}-\frac{\delta_0A}{2}\Bigr),\]
which implies the expression for $L_{\ndemi-2}$ by \eqref{relation}. 
\qed\\

\end{document}